\documentclass[a4paper, 10pt]{article}
\usepackage[colorlinks=true]{hyperref}
\hypersetup{urlcolor=blue, citecolor=red}
\usepackage[normalem]{ulem}
\usepackage{tikz}
\usepackage{fancyhdr}
\usepackage{bm}
\usepackage{cancel}
\usepackage{latexsym}
\usepackage{color}
\usepackage{amssymb,amsmath,amsfonts,bbm,pifont,upgreek,bbold,accents
}  
\usepackage[nointlimits]{esint}
\usepackage{graphicx}
\usepackage{tikz}
\newcommand\beq[1]{\begin{equation}\label{#1} }
\newcommand{\eeq}{\end{equation} }

\renewcommand{\theequation}{\arabic{section}.\arabic{equation}}

\newtheorem{theorem}{Theorem}[section]
\newtheorem{definition}{Definition}[section]
\newtheorem{proposition}{Proposition}[section]
\newtheorem{lemma}{Lemma}[section]

\newtheorem{sublemma}{Sublemma}[section]
\newtheorem{remark}{Remark}[section]
\newtheorem{notationalremark}{Notational Remark}[section]
\newtheorem{corollary}{Corollary}[section]
\newtheorem{assumption}{Assumption}[section]
\newtheorem{claim}{Claim}[section]

\newtheorem{tools}{$\negsp\negsp$}[subsection]

\newcommand\thm[1]{\begin{theorem}\label{#1}}
\newcommand\thmtwo[2]{\begin{theorem}[#1]\label{#2}}
\newcommand\ethm{\end{theorem} }
\newcommand\dfn[1]{\begin{definition}\label{#1} \rm}
\newcommand\dfntwo[2]{\begin{definition}[#1]\label{#2} \rm}
\newcommand\edfn{\end{definition} }
\newcommand\pro[1]{\begin{proposition}\label{#1}}
\newcommand\protwo[2]{\begin{proposition}[#1]\label{#2}}
\newcommand\epro{\end{proposition} }
\newcommand\lem[1]{\begin{lemma}\label{#1}}
\newcommand\lemtwo[2]{\begin{lemma}[#1]\label{#2}}
\newcommand\elem{\end{lemma} }
\newcommand\sublem[1]{\begin{sublemma}\label{#1}}
\newcommand\sublemtwo[2]{\begin{sublemma}[#1]\label{#2}}
\newcommand\esublem{\end{sublemma} }
\newcommand\rem[1]{\begin{remark}\label{#1} \rm}
\newcommand\erem{\end{remark} }
\newcommand\notrem[1]{\begin{notationalremark}\label{#1} \rm}
\newcommand\enotrem{\end{notationalremark} }
\newcommand\cor[1]{\begin{corollary}\label{#1}}
\newcommand\cortwo[2]{\begin{corollary}[#1]\label{#2}}
\newcommand\ecor{\end{corollary} }
\newcommand\asmp[1]{\begin{assumption}\label{#1}}
\newcommand\asmptwo[2]{\begin{assumption}[#1]\label{#2}}
\newcommand\easmp{\end{assumption} }
\newcommand\clm[1]{\begin{claim}\label{#1}}
\newcommand\eclm{\end{claim} }
\newcommand{\proof}{\par\medskip\noindent{\bf Proof\ }}
\expandafter\chardef\csname pre amssym.def
at\endcsname=\the\catcode`\@
\catcode`\@=11
\def\undefine#1{\let#1\undefined}
\def\newsymbol#1#2#3#4#5{\let\next@\relax
 \ifnum#2=\@ne\let\next@\msafam@\else
 \ifnum#2=\tw@\let\next@\msbfam@\fi\fi
 \mathchardef#1="#3\next@#4#5}
\def\mathhexbox@#1#2#3{\relax
 \ifmmode\mathpalette{}{\m@th\mathchar"#1#2#3}%
 \else\leavevmode\hbox{$\m@th\mathchar"#1#2#3$}\fi}
\def\hexnumber@#1{\ifcase#1 0\or 1\or 2\or 3\or 4\or 5\or 6\or 7\or
8\or
 9\or A\or B\or C\or D\or E\or F\fi}
\ifcase\@ptsize
 \font\tenmsb=msbm10
 \font\sevenmsb=msbm7
 \font\fivemsb=msbm5
\or
 \font\tenmsb=msbm10 scaled \magstephalf
 \font\sevenmsb=msbm7 scaled \magstephalf
 \font\fivemsb=msbm5  scaled \magstephalf
\or
 \font\tenmsb=msbm10 scaled \magstep1
 \font\sevenmsb=msbm7 scaled \magstep1
 \font\fivemsb=msbm5 scaled \magstep1
\fi
\newfam\msbfam
\textfont\msbfam=\tenmsb
\scriptfont\msbfam=\sevenmsb
\scriptscriptfont\msbfam=\fivemsb
\edef\msbfam@{\hexnumber@\msbfam}
\def\Bbb#1{\fam\msbfam\relax#1}
\def\widehat#1{\setboxz@h{$\m@th#1$}%
 \ifdim\wdz@>\tw@ em\mathaccent"0\msbfam@5B{#1}%
 \else\mathaccent"0362{#1}\fi}
\def\widetilde#1{\setboxz@h{$\m@th#1$}%
 \ifdim\wdz@>\tw@ em\mathaccent"0\msbfam@5D{#1}%
 \else\mathaccent"0365{#1}\fi}

\def\RIfM@{\relax\ifmmode}
\def\nonmatherr@#1{\errmessage{\string#1\space allowed only in math mode}}
\def\Bbb{\RIfM@\expandafter\Bbb@\else
 \expandafter\nonmatherr@\expandafter\Bbb\fi}
\def\Bbb@#1{{\Bbb@@{#1}}}
\def\Bbb@@#1{\fam\msbfam\relax#1}
\def\setboxz@h{\setbox\z@\hbox}
\def\wdz@{\wd\z@}
\catcode`\@=\csname pre amssym.def at\endcsname
\newcommand{\ii}{{\rm i}  }

\newcommand{\negsp}{\hspace{-.09truecm}}  

\definecolor{yellow}{rgb}{0.99, 0.93, 0.0}

\begin{document}
\title{\bf 
Proof of a conjecture by \\ H. Dullin and R. Montgomery\thanks{
 {\bf MSC2000 numbers:}
34C20, 34C25, 70F05,  70F10, 37J35, 70K43.
 {\bf Keywords:} two fixed center problem, quasi--periodic orbits, monotonicity of the periods, rotation number.
}
}
\date{}

\author{Gabriella Pinzari\footnote{Department of Mathematics, University of Padua, e--mail address: {\tt pinzari@math.unipd.it}}}
\maketitle

\begin{abstract}\footnotesize{In the framework of the planar Euler problem in the quasi--periodic regime, the formulae of the periods available in the literature are simple only on one side of their singularity. In this paper, we
 complement such formulae  with others, which result simpler on the other side. The derivation of such new formulae uses the Keplerian limit and complex analysis tools.
 As an application, we prove a conjecture by H. Dullin and R. Montgomery, which states  that such periods,  as well as their ratio, the {\it rotation number}, are monotone functions of their non--trivial first integral, at any fixed energy level. 
}
\end{abstract}

\maketitle

\tableofcontents

\newpage

\renewcommand{\theequation}{\arabic{equation}}
\setcounter{equation}{0}
\section{Introduction}\label{Purpose of the paper}

\paragraph{\bf 1.1 Overview on  the Euler problem}  The {\it two fixed centers } (or {\it Euler}--) problem is the Hamiltonian system -- firstly considered by Euler and next proven to
 be integrable by Jacobi~\cite{{Euler1957}, jacobi09} -- of the motion of a free particle undergoing gravitational attraction  by two masses in fixed positions in the Euclidean space. 
 As an integrable problem, the Euler problem has attracted the attention of many mathematicians, maybe because of the peculiarity of the process of integration, intimately connected to a spectacular regularization of collisions with {\it both} the attracting centers, which now we briefly recall (see~\cite{hiltebeitel1911, bekovO78, dullinM16, waalkensDR06} and references therein for a complete discussion).

 \noindent
 The  Hamiltonian of the Euler problem is
\begin{eqnarray}\label{oldHam}
{\rm J}=\frac{\|{\mathbf y}\|^2}{2}-\frac{{{\rm M}}}{\|{\mathbf x}+{\mathbf v}_0\|}-\frac{{\rm m}}{\|{\mathbf x}-{\mathbf v}_0\|}
\end{eqnarray}
where ${{\rm M}}$, ${\rm m}$ are the masses of the attracting bodies, ${\mathbf v}_0$ provides the direction of the mass ${\rm m}$,  which we also fix as the first axis  ${\mathbf i}$ of a reference frame: ${\mathbf v}_0=v_0\,{\mathbf i}$, with $v_0>0$;
$\|\cdot\|$ is the Euclidean norm and the gravity constant has been  put equal to $1$. Without loss of generality, we assume ${{\rm M}}\ge {\rm m}\ge 0$, but $({{\rm M}}, {\rm m})\ne (0, 0)$. Even though the result of  this paper will concern the planar case (where ``planar'' means that we fix ${\mathbf y}$, ${\mathbf x}$ to take values on a fixed plane through ${\mathbf i}$),   some formulae that we shall write hold for  the more general spatial problem, so we regard ${\mathbf y}$, ${\mathbf x}$ as objects of ${\mathbb R}^3$.  And precisely in the spatial version one observes that the Hamiltonian~\eqref{oldHam} 
 remains unvaried applying a rotation about the ${\mathbf i}$ axis, due to existence of the following first integral
 \begin{eqnarray}\label{m3}\Theta({\mathbf y}, {\mathbf x}):=x_2y_3-x_3y_2\end{eqnarray}
 consisting of the projection along the ${\mathbf i}$--axis
 of the angular momentum
\begin{eqnarray}\label{M}
{\mathbf M}({\mathbf y}, {\mathbf x}):={\mathbf x}\times {\mathbf y}
\end{eqnarray}
of the particle. Note that $\Theta$ vanishes in the planar case.\\
The solution found out by Jacobi relies on switching to new position coordinates  $\alpha$, $\beta$, $\vartheta$, with
\begin{eqnarray*}\alpha\ge 1\,,\quad |\beta|\le 1\,,\quad \vartheta\in {\mathbb T}
\end{eqnarray*}
defined as
\begin{eqnarray*}
\left\{
\begin{array}{lll}\displaystyle\alpha({\mathbf x}):=\frac{\|{\mathbf x}+{\mathbf v}_0\|+\|{\mathbf x}-{\mathbf v}_0\|}{2v_0}\\\\
\displaystyle\beta({\mathbf x}):=\frac{\|{\mathbf x}+{\mathbf v}_0\|-\|{\mathbf x}-{\mathbf v}_0\|}{2v_0}\\\\
\displaystyle\vartheta({\mathbf x}):=\arg(x_3, -x_2)
\end{array}
\right.
\end{eqnarray*}
We denote as $A({\mathbf y}, {\mathbf x})$, $B({\mathbf y}, {\mathbf x})$, $\Theta({\mathbf y}, {\mathbf x})$ are the generalized impulses conjugated to $\alpha$, $\beta$, $\vartheta$  (with $\Theta$ being precisely the quantity in~\eqref{m3}). Using such coordinates, ${\rm J}$ is carried to  the nice form (see~\cite{bekovO78} for the details)
\begin{eqnarray}\label{Jnew}{\rm J}&=&\frac{1}{v_0^2(\alpha^2-\beta^2)}\nonumber\\
&&\left(
\frac{A^2(\alpha^2-1)}{2}+\frac{B^2(1-\beta^2)}{2}+\frac{\Theta^2}{2(\alpha^2-1)}+\frac{\Theta^2}{2(1-\beta^2)}- {{\rm M}}_+v_0\alpha+ {{\rm M}}_-v_0\beta
\right)
\end{eqnarray} 
with
\begin{eqnarray}\label{M+-}{{\rm M}}_-={{\rm M}}-{\rm m}\,,\qquad {{\rm M}}_+={{\rm M}}+{\rm m}\,.\end{eqnarray}

  \vskip.1in 
  \noindent
Given the lucky  expression in~\eqref{Jnew}, Jacobi proposed to switch to the new, completely separated, Hamiltonian 
\begin{eqnarray}\label{Jreg}{\rm J}_{\rm reg}={\rm J}_+(A, \alpha,{\rm J}_0 , \Theta_0, {\rm M}_+)+{\rm J}_-(B, \beta,{\rm J}_0 , \Theta_0, {\rm M}_-)
\end{eqnarray}
where
\begin{eqnarray}\label{Jpm}
{\rm J}_+(A, \alpha,{\rm J}_0 , \Theta_0, {\rm M})
&=&\frac{A^2(\alpha^2-1)}{2}+\frac{\Theta_0^2}{2(\alpha^2-1)}-{\rm M}v_0\alpha-{\rm J}_0 v_0^2(\alpha^2-1)\nonumber\\
{\rm J}_-(B, \beta,{\rm J}_0 , \Theta_0, {\rm M})
&=&\frac{B^2(1-\beta^2)}{2}+\frac{\Theta_0^2}{2(1-\beta^2)}+{\rm M}v_0\beta-{\rm J}_0 v_0^2(1-\beta^2)
\end{eqnarray}
The motions of Hamiltonian ${\rm J}_{\rm reg}$ in~\eqref{Jreg} are related to the ones of ${\rm J}$ in~\eqref{Jnew} as follows: any orbit of ${\rm J}_{\rm reg}$ on the zero energy level in the time $\tau$ provides an orbit of ${\rm J}$ on the ${\rm J}_0$ energy level in the time $t$, related to $\tau$ via \begin{eqnarray}\label{dtaudt}\frac{d\tau}{dt}=\frac{1}{v_0^2(\alpha(t)^2-\beta(t)^2)}\,.\end{eqnarray} In addition, by the complete separability of ${\rm J}_{\rm reg}$,  the zero energy level of ${\rm J}_{\rm reg}$
enforces a further constant of motion, which we call ${\rm F}_0$, defined via
\begin{eqnarray}\label{AB}\left\{
\begin{array}{lll}
\displaystyle{\rm J}_+(A, \alpha,{\rm J}_0 , \Theta_0, {{\rm M}}_+)
=-\frac{{\rm F}_0}{2}\\\\
\displaystyle{\rm J}_-(B, \beta,{\rm J}_0 , \Theta_0, {{\rm M}}_-)
=\frac{{\rm F}_0}{2}
\end{array}
\right.
\end{eqnarray}

\begin{figure}
 \begin{center}
 \begin{tikzpicture}
\draw[densely dashed,  thick] (2,6) parabola bend (0,5)  (7.4,18.8);
\draw[white, ultra thick](2,6) parabola bend (0,5)  (2,6);
\draw[densely dashed,  thick] plot coordinates {(0,4) (2,6)};
\draw[ultra thick]  plot coordinates {(0,4) (5,0)};
\draw[densely dashed,  thick] plot coordinates {(0,4) (5,14)};
\draw[ultra thick] (4,12) parabola bend (0,8)  (4,12);
\draw[ultra thick] plot coordinates {(4,12) (7.5,19)};
\draw[ultra thick] plot coordinates {(0,4) (0,8)};
\draw[dotted, thick] plot coordinates {(0,4) (2,0)};
\draw [->] (-1,4) -- (10,4);
\draw [->] (0,0) -- (0,20);
\node at (9.5,4.3) {d};
\node at (-0.3,19.5) {f};
\node at (3.5,13) {\scriptsize${\rm f=f^+_{{\rm M}_+}({\rm d})}$};
\node at (3.3,9) {\scriptsize${\rm f=f^{\rm s}_{{\rm M}_+}({\rm d})}$};
\node at (5.0,8.5) {\scriptsize${\rm f=f^{\rm s}_{{\rm M}_-}({\rm d})}$};
\node at (5.5,0.5) {\scriptsize${\rm f=f^-_{{\rm M}_-}({\rm d})}$};
\node at (2.5,0.5) {\scriptsize${\rm f=f^-_{{\rm M}_+}({\rm d})}$};
\end{tikzpicture}
  \end{center}
  \caption{Graphical representation of the domain ${\mathbb P}$, in the plane ${\rm (d, f)}$. The singular lines $\rm f=f^{\rm s}_{{\rm M}_\pm}({\rm d})$ (dashed) and the boundary lines of ${\mathbb P}$ (thick) are reported. The  lower boundary $\rm f=f^{-}_{{\rm M}_-}({\rm d})$ of $\mathbb P$ is also a boundary line of ${\mathbb Q}_{\rm M_-}$, while 
  the line
  $\rm f=f^{-}_{{\rm M}_+}({\rm d})$ (dotted), lower boundary line of ${\mathbb Q}_{\rm M_+}$, is external to $\mathbb P$. It
  has been reported for comparison with Fig.~\ref{domains QM}. }\label{domainsOLD}
 \end{figure}
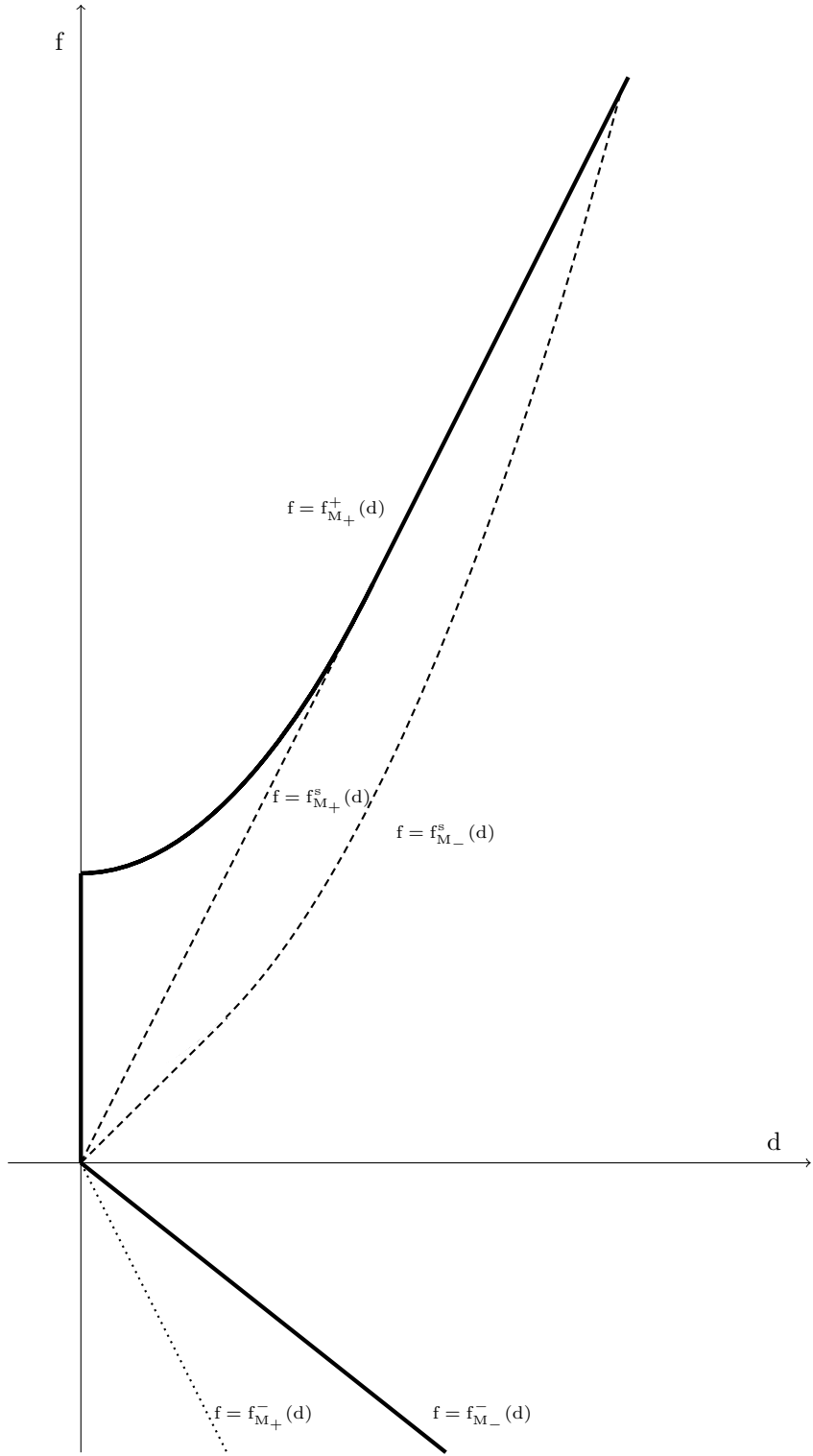

\paragraph {\bf 1.2 Problem and result} 
In this paper, we are interested  to the periods (hereafter called {\it Jacobi periods}, and denoted as $\tau_+$, $\tau_-$) of the Euler problem, in the quasi--periodic regime.
The classical analysis that we have reviewed in the previous section is fit to completely describe the set ${\mathbb P}$ of parameters $({\rm J}_0, {\rm F}_0, \Theta_0)$ corresponding to periodic orbits, and provide an expression of their periods, in form of elliptic integrals. In this paper, we aim to complement the formulae for $\tau_+$, $\tau_-$ existing in the literature with others, so as to write down the simplest expressions possible. For simplicity, we focus on the planar case ($\Theta_0=0$)  (even though we do not see obstructions to treat the general case with the same techniques here) and switch from ${\rm J}_0$ and ${\rm F}_0$ to the quantities
 \begin{eqnarray}\label{new parameters}
{{\rm d}}:=-4 v_0 {\rm J}_0\,,\qquad {{\rm f}}:=-2{\rm J}_0{\rm F}_0
\end{eqnarray}
Observe that,
by the assumptions on $\rm m$, $\rm M$, the mass parameters ${\rm M}_-$, ${\rm M}_+$ in~\eqref{M+-} verify
  \begin{eqnarray}\label{M+M+}
  0\le {\rm M}_-\le {\rm M}_+\,, \quad ({\rm M}_-,  {\rm M}_+)\ne (0, 0)
  \end{eqnarray}
with ${\rm M}_-=0$ corresponding to ${\rm m}={\rm M}$, while   ${\rm M}_-={\rm M}_+$ to ${\rm m}=0$. From now on, we shall implicitly assume~\eqref{M+M+}, without further mention.
Then the set ${\mathbb P}$ is a convex  set, twice disconnected, having the form (see Fig.~\ref{domainsOLD})
\begin{eqnarray}\label{mathbbP}
{\mathbb P}=\left\{({{\rm d}}, {{\rm f}})\in {\mathbb R}^2:\ 
{{\rm d}}>0\,,\quad {{\rm f}}^-_{{\rm M}_-}({{\rm d}})<{{\rm f}}< {{\rm f}}^+_{{\rm M}_+}({{\rm d}})\,,\ {{\rm f}}\ne {{\rm f}}^{\rm s}_{{\rm M}_-}({{\rm d}})\,,\ {{\rm f}}\ne {{\rm f}}^{\rm s}_{{\rm M}_+}({{\rm d}})
\right\}
\end{eqnarray}
 \!with suitable ${{\rm f}}^-_{{\rm M}}({{\rm d}})$, $ {{{\rm f}}}^{\rm s}_{{\rm M}}({{\rm d}})$, ${{\rm f}}^+_{{\rm M}}({{\rm d}})$
precisely defined in the course of the paper (see Equation~\eqref{FmaxFsing} below)
and verifying
\begin{eqnarray}\label{fpms}{{\rm f}}^-_{{\rm M}_+}({{\rm d}})\le{{\rm f}}^-_{{\rm M}_-}({{\rm d}})\le  {{{\rm f}}}^{\rm s}_{{\rm M}_-}({{\rm d}})\le  {{{\rm f}}}^{\rm s}_{{\rm M}_+}({{\rm d}})\le {{\rm f}}^+_{{\rm M}_+}({{\rm d}})\quad \forall\ {\rm d}>0\,.\end{eqnarray}
We now define, for each ${\rm M}\in \{{\rm M}_-, {\rm M}_+\}$, the set
 \begin{eqnarray}\label{QMNEW}{{\mathbb Q}}_{\rm M}:=\left\{({{\rm d}}, {{\rm f}})\in {\mathbb R}^2:\ 
{{\rm d}}>0\,,\quad {{\rm f}}>{{\rm f}}^-_{{\rm M}}({{\rm d}})\,,\ {\rm f}\ne{{\rm f}}^{\rm s}_{{\rm M}}({{\rm d}})
\right\}\end{eqnarray}
   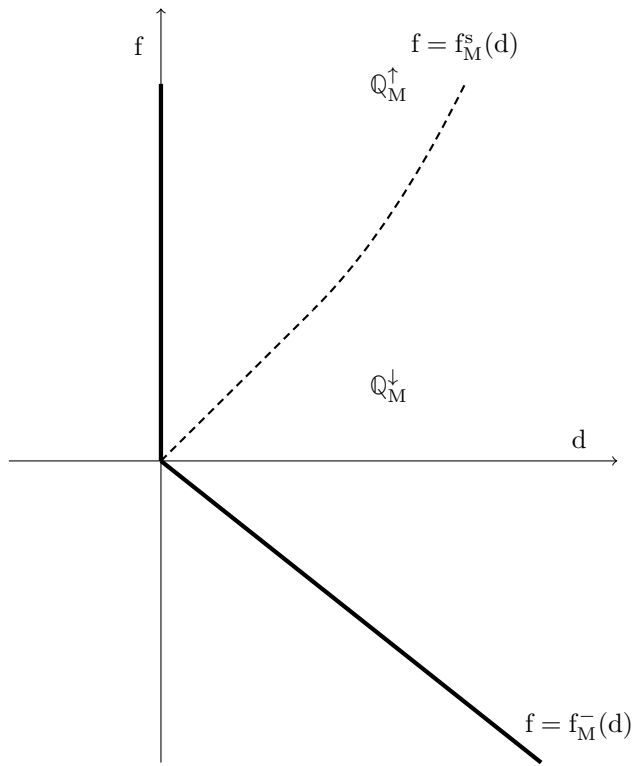
\begin{figure}
 \begin{center}
 \begin{tikzpicture}
\draw[densely dashed, thick] (2,6) parabola bend (0,5)  (4,9);
\draw[white, ultra thick] (2,6) parabola bend (0,5)  (2,6);
\draw[densely dashed,  thick] plot coordinates {(0,4) (2,6)};
\draw[ultra thick]  plot coordinates {(0,4) (5,0)};
\draw [->] (-2,4) -- (6,4);
\draw [->] (0,0) -- (0,10);
\draw [ultra thick] (0,4) -- (0,9);
\node at (5.5,4.3) {d};
\node at (-0.3,9.5) {f};
\node at (4.0,9.5) {${\rm f=f^{\rm s}_{\rm M}({\rm d})}$};
\node at (5.5,0.5) {${\rm f=f^-_{\rm M}({\rm d})}$};
\node at (3,5) {${\mathbb Q^{\downarrow}_{\rm M}}$};
\node at (3,9) {${\mathbb Q^{\uparrow}_{\rm M}}$};
\end{tikzpicture}
  \end{center}
  \caption{Graphical representation of the domain ${{\mathbb Q}}_{\rm M}$, and of its sub--domains ${{\mathbb Q}}^{\downarrow}_{\rm M}$, ${{\mathbb Q}}^{\uparrow}_{\rm M}$, in the plane ${\rm (d, f)}$. The ``singular line'' ${\rm f=f^{\rm s}_{\rm M}({\rm d})}$ (dashed) and the 
  ``minimum line'' ${\rm f=f^-_{\rm M}({\rm d})}$ are reported.  The domain $\mathbb P$ in Fig.~\ref{domainsOLD} is a subset of ${{\mathbb Q}}^{\uparrow}_{\rm M_+}\cap {{\mathbb Q}}^{\uparrow}_{\rm M_+}$.}\label{domains QM}
 \end{figure}
\!By~\eqref{mathbbP} and~\eqref{fpms}, we have: ${{\mathbb Q}}_{{\rm M}_+}\,,\ {{\mathbb Q}}_{{\rm M}_-}\supset {\mathbb P}$.
We split each ${{\mathbb Q}}_{\rm M}$ as
${{\mathbb Q}}_{\rm M}={{\mathbb Q}}_{\rm M}^\uparrow\cup {{\mathbb Q}}_{\rm M}^\downarrow$,
where ${{\mathbb Q}}_{\rm M}^\uparrow$ is the part of ${{\mathbb Q}}_{\rm M}$ above the line $\{{\rm f}={\rm f}^{\rm s}_{\rm M}({\rm d})\}$, while ${{\mathbb Q}}_{\rm M}^\downarrow$ is the one below (see Fig.~\ref{domains QM}). On each part of ${{\mathbb Q}}_{\rm M}$, we define the function

\begin{eqnarray}\label{thetaMupanddown}
\theta_{\rm M}({{\rm d}}, {{\rm f}}):=\left\{\begin{array}{lll} \displaystyle \theta^{\uparrow}_{\rm M}:=2\sqrt{\frac{2\rm d}{v_0}}\int_{-1}^{1}\frac{dx}{ \sqrt{(1-x^2)\big({\rm d}^2x^2-4{\rm d}{{\rm M}}x+4{\rm f}-{\rm d}^2\big)}}\quad &{\rm if}\ ({{\rm d}}, {{\rm f}})\in
{{\mathbb Q}}^{\uparrow}_{\rm M}\\\\
\displaystyle \theta^{\downarrow}_{\rm M}:=\sqrt{\frac{2\rm d}{v_0}}\int_{0}^{2}\frac{dz}{\sqrt{z(2-z)({\rm M}^2 z^2-2{\rm f} z+{\rm d}^2)}}
\quad &{\rm if}\ ({{\rm d}}, {{\rm f}})\in {{\mathbb Q}}^\downarrow_{\rm M}
\end{array}
\right.
\end{eqnarray}

\noindent
In this paper we prove that
\begin{theorem}\label{extendperiods}  
The function $\theta_{\rm M}$ is well--defined on ${{\mathbb Q}}_{\rm M}$
and, moreover, the following identity holds
\begin{eqnarray*}
\tau_\pm={\theta}_{{\rm M}_\pm}|_{\mathbb P}\quad \forall\ ({\rm d}, {\rm f})\in {\mathbb P}\,.
\end{eqnarray*}
\end{theorem}
We remark that the well posedness of $\theta^\uparrow_{{\rm M}}$ on ${{\mathbb Q}}^\uparrow_{\rm M}$ and the validity
of  
\begin{eqnarray*}
\tau_\pm={\theta}^\uparrow_{{\rm M}_\pm}|_{\mathbb P}\quad  \forall\ ({\rm d}, {\rm f})\in {\mathbb P}\cap {{\mathbb Q}}^\uparrow_{\rm M_\pm}
\end{eqnarray*}
 is a mere consequence of the classical analysis, as shown in  Section~\ref{App: classical periods}.  The novelty of the paper relies on  the well posedness of $\theta^\downarrow_{{\rm M}}$ on ${{\mathbb Q}}^\downarrow_{\rm M}$ and the validity
of  
 \begin{eqnarray*}
\tau_\pm={\theta}^\downarrow_{{\rm M}_\pm}|_{\mathbb P}\quad  \forall\ ({\rm d}, {\rm f})\in {\mathbb P}\cap {{\mathbb Q}}^\downarrow_{\rm M_\pm}\,.
\end{eqnarray*}
This is discussed in Sections~\ref{Proof of Theorem1.2} and~\ref{A simpler formula}.

\paragraph{1.3 An application} As an application of Theorem~\ref{extendperiods}, we shall prove a conjecture posed\footnote{Precisely, the conjecture regards the monotonicity of ${\rm f}\to \tau_{+}$ and of the functions
${\rm f}\to\omega_{S}$, ${\rm f}\to\omega_{P}$ defined below. The remaining cases have been 
 proven in~\cite{dullinM16}, using the theory of elliptic functions.  } by H. Dullin and R. Montgomery  in~\cite{dullinM16}. Besides the functions $\tau_-({{\rm d}}, {{\rm f}})$, $\tau_+({{\rm d}}, {{\rm f}})$, consider also  their ratio, the {\it rotation number}:
  \begin{eqnarray*}\omega({{\rm d}}, {{\rm f}}):=\frac{\tau_-({{\rm d}}, {{\rm f}})}{\tau_+({{\rm d}}, {{\rm f}})}\,.\end{eqnarray*}
We focus on the dependence of the functions $\tau_-({{\rm d}}, {{\rm f}})$, $\tau_+({{\rm d}}, {{\rm f}})$ and $\omega({{\rm d}}, {{\rm f}})$ on $\rm f$. To this end, given ${\mathbb A}\subset {\mathbb R}^2$, we denote as ${{\mathbb A}}_{\rm d}$ the set of $\rm f$ such that $({\rm d}, {\rm f})\in {{\mathbb A}}$, for a fixed $\rm d>0$, and observe that 
${\mathbb P}_{\rm d}$ is twice disconnected along ${\rm f}={\rm f}^{\rm s}_{\rm M_-}({\rm d})$, ${\rm f}={\rm f}^{\rm s}_{\rm M_+}({\rm d})$ when $0\le \rm M_-<\rm M_+$; once when $\rm M_-=\rm M_+$.
We  prove that
\begin{theorem}\label{main2}~The functions $\tau_-({{\rm d}}, \cdot)$, $\tau_+({{\rm d}}, \cdot)$ and $\omega({{\rm d}}, \cdot)$ are differentiable and monotone on each connected component of ${\mathbb P}_{\rm d}$.

\end{theorem}

\noindent
Besides proving a problem left open in~\cite{dullinM16}, 
Theorem~\ref{main2} has been recently applied in the framework of Rabinowitz Floer homology~\cite{frauenfelder2023, TakeuchiZ2023}. In addition, we foresee applications to close-to-be integrable systems having the two-center problem as limiting case, like the  problem of three or more bodies. Various results in this directions have been obtained by the author and collaborators: see~\cite{pinzari20b} for the description of a possible setting;~ \cite{pinzari20a, chenP2021} for rigorous results;~\cite{diruzzaDP20, diruzzaP2022} for numerical studies. 
Moreover, as the monotonicity of the periods is somewhat  a convexity assertion, we  foresee connections to  Nekhorossev theory~\cite{nehorosev77, nehorosev79, poschel93, niederman04, guzzoCB16}, while the monotonicity of the rotation number  should allow interactions with 
weak {\sc kam} theory: see~\cite{mather1985, mather1991, moeckel2015, bonannoM2022, siconolfiS2023} for general information;~\cite{chang2022} for an application of Aubry--Mather theory to the Euler problem. All such possible directions of research are here mentioned just as possible hints, being definitely far from the purposes of the paper. Here we  limit to show how Theorem~\ref{main2} follows from Theorem~\ref{extendperiods}.
\proof {\bf of Theorem~\ref{main2}}
The differentiability of $\tau_-(\rm d, \cdot)$, $\tau_+(\rm d, \cdot)$ and $\omega(\rm d, \cdot)$ and monotonicity of $\tau_-(\rm d, \cdot)$ and $\tau_+(\rm d, \cdot)$ on ${\mathbb P}_{\rm d}$ are an immediate consequence\footnote{As a matter of fact,  when $\rm M_-\ne \rm M_+$, Theorem~\ref{extendperiods} and the formula in~\eqref{thetaMupanddown} show a little more, namely, that  $\tau_-(\rm d, \cdot)$ and $\tau_+(\rm d, \cdot)$  are respectively differentiable and monotone for  ${\rm f}\in{\mathbb P}_{\rm d}\cup\{{\rm f}={\rm f}^{\rm s}_{\rm M_+}({\rm d})\}$, ${\rm f}\in{\mathbb P}_{\rm d}\cup\{{\rm f}={\rm f}^{\rm s}_{\rm M_-}({\rm d})\}$, respectivey.}  of Theorem~\ref{extendperiods} and of the formula in~\eqref{thetaMupanddown}. 
 Concerning $\omega(\rm d, \cdot)$ for ${\rm f}\in{\mathbb P}_{\rm d}$, it is sufficient consider the case $0\le {{\rm M}}_-<{{\rm M}}_+$, as when ${{\rm M}}_-={{\rm M}}_+$,  we are in the case of the Kepler problem, so $\tau_+=\tau_-$, whence $\omega({{\rm d}}, {{\rm f}})\equiv 1$ is trivially monotone.  As
${{\mathbb Q}}^\downarrow_{{{\rm M}}_-}\subset{{\mathbb Q}}^\downarrow_{{{\rm M}}_+}$, $ {{\mathbb Q}}^\uparrow_{{{\rm M}}_+}\subset{{\mathbb Q}}^\uparrow_{{{\rm M}}_-}$ for all $0\le {{\rm M}}_-<{{\rm M}}_+$,
then $\omega(\rm d, \rm f)$ is given by

\begin{eqnarray*}
\omega({{\rm d}}, {{\rm f}})=\left\{\begin{array}{lll}
\displaystyle \omega_S({{\rm d}}, {{\rm f}}):=\frac{\theta^\downarrow_{{{\rm M}}_-}({{\rm d}}, {{\rm f}})}{\theta^\downarrow_{{{\rm M}}_+}({{\rm d}}, {{\rm f}})}\quad &{\rm if}\ ({{\rm d}}, {{\rm f}})\in {{\mathbb D}}_S:={{\mathbb Q}}^\downarrow_{{{\rm M}}_-}\\
\displaystyle \omega_L({{\rm d}}, {{\rm f}}):=\frac{\theta^{\uparrow}_{{{\rm M}}_-}({{\rm d}}, {{\rm f}})}{\theta^\downarrow_{{{\rm M}}_+}({{\rm d}}, {{\rm f}})}\quad &{\rm if}\ ({{\rm d}}, {{\rm f}})\in {{\mathbb D}}_L:={{\mathbb Q}}^\uparrow_{{{\rm M}}_-}\cap {{\mathbb Q}}^\downarrow_{{{\rm M}}_+}\\
\displaystyle \omega_P({{\rm d}}, {{\rm f}}):=\frac{\theta^{\uparrow}_{{{\rm M}}_-}({{\rm d}}, {{\rm f}})}{\theta^{\uparrow}_{{{\rm M}}_+}({{\rm d}}, {{\rm f}})}\quad &{\rm if}\  ({{\rm d}}, {{\rm f}})\in {{\mathbb D}}_P:={{\mathbb Q}}^\uparrow_{{{\rm M}}_+}
\end{array}
\right.
\end{eqnarray*}
where, for ease of comparison, we have used the same notations as in~\cite{dullinM16}, for what concerns the subscripts $S$, $L$, $P$.
We shall prove the following assertions which, incidentally, correspond to the second figure of the ones numbered as 12 in~\cite{dullinM16}, but  obtained numerically therein.

\begin{itemize}
\item[\rm(S)] The function $ \omega_S({{\rm d}}, \cdot)$ increases from a positive finite value to $+\infty$ (as ${\rm f}\to {\rm f}^{\rm s}_{\rm M_-}(\rm d)^-$)
while ${{\rm f}}\in {{\mathbb D}}_{S, \rm d}$.
\item[\rm(L)] The function $ \omega_L({{\rm d}}, \cdot)$ decreases from $+\infty$ (as ${\rm f}\to {\rm f}^{\rm s}_{\rm M_-}(\rm d)^+$) to $0$ (as ${\rm f}\to {\rm f}^{\rm s}_{\rm M_+}(\rm d)^-$)
while ${{\rm f}}\in {{\mathbb D}}_{S, \rm d}$.
\item[\rm (P)]  The function $\omega_P({{\rm d}}, {{\rm f}})$ increases from $0$ (as ${\rm f}\to {\rm f}^{\rm s}_{\rm M_+}(\rm d)^+$) to a positive finite value
while ${{\rm f}}\in {{\mathbb D}}_{P, \rm d}$.
\end{itemize}
\!The only assertions among the one listed above which deserve a discussion concern the  increasing monotonicity of $\omega_S(\rm d, \cdot)$ and of $\omega_P(\rm d, \cdot)$, as
the limiting values of $\omega_S(\rm d, {\rm f})$, $\omega_L(\rm d, {\rm f})$ and $\omega_P(\rm d, {\rm f})$ as ${\rm f}\to {\rm f}^{\rm s}_{\rm M_\pm}(\rm d)^\pm$  as well as the decreasing monotonicity of $\omega_L(\rm d, \cdot)$ are a trivial consequence of the formula~\eqref{thetaMupanddown} and Theorem~\ref{extendperiods}. Let us check the monotonicity of $\omega_S(\rm d, \cdot)$; the case $\omega_P(\rm d, \cdot)$ being analogous. 
From the formula  \begin{eqnarray*}
\partial_{\rm f}\omega_S({{\rm d}}, {{\rm f}})=\partial_{\rm f}\left(\frac{\theta^\downarrow_{{{\rm M}}_-}({{\rm d}}, {{\rm f}})}{\theta^\downarrow_{{{\rm M}}_+}({{\rm d}}, {{\rm f}})}\right)=\frac{{\theta^\downarrow_{{{\rm M}}_+}}
{\partial_{\rm f}\theta^\downarrow_{{{\rm M}}_-}}-{\theta^\downarrow_{{{\rm M}}_-}}{\partial_{\rm f}\theta^\downarrow_{{{\rm M}}_+}}}{
\theta^\downarrow_{{{\rm M}}_+}({{\rm d}}, {{\rm f}}
\big)^2}
\end{eqnarray*}
it comes that we need to prove
\begin{eqnarray*}{\theta^\downarrow_{{{\rm M}}_+}}
{\partial_{\rm f}\theta^\downarrow_{{{\rm M}}_-}}-{\theta^\downarrow_{{{\rm M}}_-}}{\partial_{\rm f}\theta^\downarrow_{{{\rm M}}_+}}>0\,.\end{eqnarray*}
Using~\eqref{thetaMupanddown}, we have that the left hand side of this inequality is given by
\begin{eqnarray*}
\iint_{[0, 2]^2}&&\frac{x}{\sqrt{x(2-x)({\rm M}_-^2x^2-2{{\rm f}}x+{{\rm d}}^2)}\sqrt{y(2-y)({\rm M}_+^2y^2-2{{\rm f}}y+{{\rm d}}^2)}}\nonumber\\
&&\left(\frac{1}{{\rm M}_-^2x^2-2{{\rm f}}x+{{\rm d}}}-\frac{1}{{\rm M}_+^2x^2-2{{\rm f}}x+{{\rm d}}^2}\right)dxdy
\end{eqnarray*}
The function under the integral is manifestly positive for all $0\le{\rm M}_-<{\rm M}_+$, hence so il the integral itself.
$\quad \square$
\noindent

\paragraph{1.4 On the proof of Theorem~\ref{extendperiods}}
We shall prove a more general result (Proposition~\ref{prop: new period} below) which, in combination with the classical representation of the periods (Proposition~\ref{thm: classical periods}) implies, in particular,  Theorem~\ref{extendperiods}. In turn, the proof of Proposition~\ref{prop: new period} will be developed in two steps, corresponding to Sections~\ref{Proof of Theorem1.2} and~\ref{A simpler formula}, respectively. Such two steps correspond to the proof of Proposition~\ref{prop: new period}  for values of the parameters $({\rm d}, {\rm f})$ in~\eqref{new parameters} belonging to two different subsets of ${{\mathbb Q}}_{\rm M}$ (subsets whose union gives ${{\mathbb Q}}_{\rm M}$, possibly deprived of a finite number of curves). Such two proofs are completely different one from the other, as the proof 
in Sections~\ref{Proof of Theorem1.2} uses dynamical arguments, while the one in Section 
\ref{A simpler formula} is purely analytic, and uses  the theory of elliptic functions. Here we provide a brief account of the former only.  All starts with observing  that, due to the separability of ${\rm J}_{\rm reg}$ in~\eqref{Jreg},
  $\tau_+$ and $\tau_-$ depend on the masses ${\rm m}$ and ${{\rm M}}$ only via the combinations ${{\rm M}}_+$, ${{\rm M}}_-$ in~\eqref{M+-}. 
 We write $\tau_+=\tau_+({\rm M_+})$, $\tau_-=\tau_-({\rm M_-})$. Then it is quite natural  to try to reconstruct $\tau_+$, $\tau_-$   from the expressions of the same quantities of the case when ${\rm m}=0$ and ${{\rm M}}\in \{{{\rm M}}_+$, ${{\rm M}}_-\}$.
But 
 (provided that ${\rm M}_-\ne0$) this is nothing else that looking at the orbits of the Kepler Hamiltonian ${\rm K}$ (given in~\eqref{kepler} below) with sun mass equal to ${\rm M}_-$, ${\rm M}_+$, respectively,
on the level set defined by ${\rm F}_0$ and  ${\rm K}={\rm J}_0$, along the time $\tau$ in~\eqref{dtaudt}. As quasi--periodic orbits of the Kepler problem are actually periodic, 
the values of $\tau_+({\rm M}_+)$ and $\tau_-({\rm M}_-)$ for such case
 coincide, and, calling $\tau^{\rm Kep}_{\rm M}$ such common value, it turns out that the domains
 of $\tau^{\rm Kep}_{\rm M_+}$, $\tau^{\rm Kep}_{\rm M_-}$ are
are strictly smaller than the original domain of $\tau_+({\rm M}_+)$, $\tau_-({\rm M}_-)$. 
Were it not enough, we shall be able to derive a really simple formula  for $\tau^{\rm Kep}_{\rm M}$ (corresponding to the second line in~\eqref{thetaMupanddown}) only on a smaller subset of the domain of $\tau^{\rm Kep}_{\rm M}$. 
To extend the formulae of $\tau^{\rm Kep}_{\rm M_+}$, $\tau^{\rm Kep}_{\rm M_+}$ on the the original domains of $\tau_+({\rm M}_+)$, $\tau_-({\rm M}_-)$, and show that such extended formulae coincide with the correct formulae of $\tau_+({\rm M_+})$, $\tau_-({\rm M_-})$, we shall need an analytic discussion, which we provide in Section~\ref{A simpler formula}.

\section{Classical formulae for Jacobi periods}\label{App: classical periods} 
Define the function
\begin{eqnarray}\label{TupFORMULA}
\tau_{\rm M}:=\left\{\begin{array}{lll} \displaystyle \tau^{\uparrow}_{\rm M}:=2\sqrt{\frac{2\rm d}{v_0}}\int_{-1}^{1}\frac{dx}{ \sqrt{(1-x^2)\big({\rm d}^2x^2-4{\rm d}{{\rm M}}x+4{\rm f}-{\rm d}^2\big)}}\quad &{\rm if}\ ({{\rm d}}, {{\rm f}})\in
{{\mathbb Q}}^{\uparrow}_{\rm M}\\\\
\displaystyle \tau^{\downarrow}_{\rm M}:=2\sqrt{\frac{2\rm d}{v_0}}\int_{-1}^{x_{\rm M}({\rm d}, {\rm f})}\frac{dx}{ \sqrt{(1-x^2)\big({\rm d}^2x^2-4{\rm d}{{\rm M}}x+4{\rm f}-{\rm d}^2\big)}}\ &{\rm if}\ ({{\rm d}}, {{\rm f}})\in {{\mathbb Q}}^\downarrow_{\rm M}
\end{array}
\right.
\end{eqnarray}
where $x_{\rm M}({\rm d}, {\rm f})$ is the minimum root of ${\rm d}^2x^2-4{\rm d}{{\rm M}}x+4{\rm f}-{\rm d}^2$
and ${{\mathbb Q}}^\downarrow_{\rm M}$, ${{\mathbb Q}}^\uparrow_{\rm M}$ are defined as in the introduction.
The purpose of this paper is to show that the  
 classical analysis reviewed in Section~\ref{Purpose of the paper} leads to the following 
\begin{proposition}\label{thm: classical periods}
The function $\tau_{\rm M}$ is well--defined on ${{\mathbb Q}}_{\rm M}$, tends to $+\infty$ as ${\rm f}\to {\rm f}^{\rm s}_{\rm M}({\rm d})^\pm$ and, moreover, the following identity holds
\begin{eqnarray*}
\tau_\pm=\tau_{{\rm M}_\pm}|_{\mathbb P}\quad \forall\ ({\rm d}, {\rm f})\in {\mathbb P}
\end{eqnarray*}
\end{proposition}
\begin{remark}\rm
Remark that $\tau_{\rm M}^\uparrow=\theta_{\rm M}^\uparrow$, while $\tau_{\rm M}^\downarrow$ is much more complicated compared to $\theta_{\rm M}^\downarrow$ (compare~\eqref{thetaMupanddown} and~\eqref{TupFORMULA}). In particular, the proof of Theorem~\ref{main2}  based on the formulae in~\eqref{TupFORMULA} does not seem simple.
\end{remark}

\vskip.1in
\noindent
The separability of the Hamiltonian ${\rm J}_{\rm reg}$ in~\eqref{Jreg} allows to consider the motions of the variables $\alpha $ and $\beta$ independently one of the other. In particular, one can first define the two respective
``Hill sets'' 
{\small
\begin{eqnarray}\label{Hill}
{{{\mathbb H}}}_+:=\Big\{
({\rm J}_0, {\rm F}_0, \Theta_0)\in{\mathbb R}^3 :\ {{\mathbb A}}({\rm J}_0, {\rm F}_0, \Theta_0)\ne \emptyset\Big\}\,,\ 
{{{\mathbb H}}}_-:=\Big\{
({\rm J}_0, {\rm F}_0, \Theta_0)\in{\mathbb R}^3 :\ {{\mathbb B}}({\rm J}_0, {\rm F}_0, \Theta_0)\ne \emptyset\Big\}\nonumber\\
\end{eqnarray}}
\!\!\!where
\begin{eqnarray*}
&&{{\mathbb A}}({\rm J}_0, {\rm F}_0, \Theta_0):=\Big\{\alpha\in (1, +\infty):\ (\alpha^2-1)\big(2{\rm J}_0 v_0^2\alpha^2+2 {{\rm M}}_+v_0\alpha-(2{\rm J}_0 v_0^2+{\rm F}_0)\big)-\Theta_0^2\ge 0\Big\}\nonumber\\
&&{{\mathbb B}}({\rm J}_0, {\rm F}_0, \Theta_0):=\Big\{
\beta\in (-1, 1):\ (1-\beta^2)\big(-2{\rm J}_0 v_0^2\beta^2-2 {{\rm M}}_-v_0\beta+(2{\rm J}_0 v_0^2+{\rm F}_0)\big)-\Theta_0^2\ge 0
\Big\}
\end{eqnarray*}
and next consider the ``Hill set of ${\rm J}_{\rm reg}$'', defined as the intersection
\begin{eqnarray*}{{\mathbb H}}:={{{\mathbb H}}}_+\cap {{{\mathbb H}}}_-\,.\end{eqnarray*}
Likewise, one can depict the set of parameters giving rise to periodic orbits for ${\rm J}_{\rm reg}$ as the intersection of analogous sets for the two planes, separately:
\begin{eqnarray}\label{P0-OLD}
{{\mathbb P}}({\rm M}_+, {\rm M}_-):={{\mathbb P}}_+({\rm M}_+)\cap {{\mathbb P}}_-({\rm M}_-)
\end{eqnarray}
where 
\begin{eqnarray*}
{{\mathbb P}}_{+}({\rm M}_+)&:=&\Big\{
({\rm J}_0, {\rm F}_0, \Theta_0)\in{{{\mathbb H}}}_+ :\ \rm the\ \alpha\textrm{--}level\ curves\ are\ smooth, \ closed\ and\ connected
\Big\}\nonumber\\
{{\mathbb P}}_{-}({\rm M}_-)&:=&\Big\{
({\rm J}_0, {\rm F}_0, \Theta_0)\in{{{\mathbb H}}}_-:\ \rm the\ \beta\textrm{--}level\ curves\ are\ smooth, \ closed\ and\ connected
\Big\}
\end{eqnarray*}
where ``$\alpha$, $\beta$ {\it level curves}'' are defined as 
the curves in the planes $(\alpha, \alpha')$, $(\beta, \beta')$ defined\footnote{Combining the evolution equation for $\alpha'$ obtained from the Hamilton equations of ${\rm J}_{+}$ and fixing the level~\eqref{AB} one obtains
\begin{eqnarray*}
\alpha'&=&A(\alpha^2-1)=\pm \sqrt{(\alpha^2-1)\big(2{\rm J}_0 v_0^2\alpha^2+2 {{\rm M}}_+v_0\alpha-(2{\rm J}_0 v_0^2+{\rm F}_0)\big)-\Theta_0^2}
\end{eqnarray*}
The equation for $\beta'$ is obtained similarly.
} via 
\begin{eqnarray}\label{eq for alpha}
\alpha'&=&\pm \sqrt{(\alpha^2-1)\big(2{\rm J}_0 v_0^2\alpha^2+2 {{\rm M}}_+v_0\alpha-(2{\rm J}_0 v_0^2+{\rm F}_0)\big)-\Theta_0^2}\\
\label{eq for beta}\beta'&=&\pm \sqrt{(1-\beta^2)\big(-2{\rm J}_0 v_0^2\beta^2-2 {{\rm M}}_-v_0\beta+(2{\rm J}_0 v_0^2+{\rm F}_0)\big)-\Theta_0^2}
\end{eqnarray}
for $({\rm J}_0, {\rm F}_0, \Theta_0)\in {\mathbb H}_+$, ${\mathbb H}_-$, respectively. 
The classical analysis reviewed in Section~\ref{Purpose of the paper} leads to the following expressions for the periods
\begin{eqnarray}\label{TaTb}
\tau_+({\rm J}_0, {\rm F}_0, \Theta_0)&=&2\int_{\alpha_{\rm min}({\rm J}_0, {\rm F}_0, \Theta_0)}^{\alpha_{\rm max}({\rm J}_0, {\rm F}_0, \Theta_0)}\nonumber\\
&&\frac{d\alpha}{ \sqrt{(\alpha^2-1)\big(2{\rm J}_0 v_0^2\alpha^2+2{{\rm M}}_+v_0\alpha-({\rm F}o+2{\rm J}_0 v_0^2)\big)-\Theta_0^2}}\nonumber\\\nonumber\\
\tau_-({\rm J}_0, {\rm F}_0, \Theta_0)&=&2\int_{\beta_{\rm min}({\rm J}_0, {\rm F}_0, \Theta_0)}^{\beta_{\rm max}({\rm J}_0, {\rm F}_0, \Theta_0)}\nonumber\\
 &&\frac{d\beta}{ \sqrt{(1-\beta^2)\big(-2{\rm J}_0 v_0^2\beta^2-2{{\rm M}}_-v_0\beta+({\rm F}o+2{\rm J}_0 v_0^2)\big)-\Theta_0^2}}
 \end{eqnarray}
where $\alpha_{\rm min}$, $\alpha_{\rm max}$ ($\beta_{\rm min}$, $\beta_{\rm max}$) are the intersections of the $\alpha$ ($\beta$) smooth level curves with the axis $\alpha'=0$ ($\beta'=0$).
From now on, we focus on the case $\Theta_0=0$ and use the parameters $({{\rm d}}, {{\rm f}})$ in~\eqref{new parameters}, instead of ${\rm J}_0$, ${\rm F}_0$.    Let ${{\rm f}}^-_{{\rm M}}({{\rm d}})$, ${{\rm f}}^+_{{\rm M}}({{\rm d}})$ and ${{\rm f}}^{\rm s}_{{\rm M}}({{\rm d}})$ be defined as 
\begin{eqnarray}\label{FmaxFsing}
 {{\rm f}}^-_{{\rm M}}({{\rm d}})&:=&-{\rm M}{\rm d}\nonumber\\
 {{{\rm f}}}^{\rm s}_{{\rm M}}({{\rm d}})
&:=&
\left\{
\begin{array}{lll}
 {{\rm M}}{{\rm d}}\quad &{\rm if}\quad 0<{{\rm d}} \le 2 {\rm M}\\
 {{\rm M}}^2+\frac{{{\rm d}}^2}{4}
 &{\rm if}\quad {{\rm d}}
 > 2 {\rm M}
\end{array}
\right.\ {{{\rm f}}}^{+}_{{\rm M}}({{\rm d}})
:=
\left\{
\begin{array}{lll}
 {{\rm M}}{{\rm d}}
\quad &{\rm if}\quad {{\rm d}} > 2 {\rm M}\\
 {{\rm M}}^2+\frac{{{\rm d}}^2}{4}
 &{\rm if}\quad 0<{{\rm d}} \le 2 {\rm M}\end{array}
\right.\end{eqnarray}
Here, we have assigned the superscripts  ``+' ,  ``-'' and ``s'' to recall the words ``maximum'' , ``minimum'' and ``singular'', as in the Keplerian case (${\rm M}_+={\rm M}_-={\rm M}$), ${{{\rm f}}}^{+}_{{\rm M}}({{\rm d}})$, ${{{\rm f}}}^{-}_{{\rm M}}({{\rm d}})$ are, respectively, the maximum and the minimum possible values of $\rm f$, while  ${{{\rm f}}}^{\rm s}_{{\rm M}}({{\rm d}})$ is the unique value of $\rm f$ such that the Jacobi period is not defined (compare Proposition~\ref{Keplerian periods} below).
\begin{proposition}\label{classical Jacobi periods}
The sets{ ${\mathbb P}_+({\rm M}_+)$, ${\mathbb P}_-({\rm M}_-)$ are given by}
 \begin{eqnarray}\label{P0-}{\mathbb P}_{+}({\rm M}_+)&=&\left\{({{\rm d}}, {{\rm f}})\in {\mathbb R}^2:\ {{\rm d}}\ge 0\,,\qquad  {{\rm f}}
\le {{\rm f}}^+_{{{\rm M}}_+}({{\rm d}})\,,\qquad
{{\rm f}}
\ne {{\rm f}}^{\rm s}_{{{\rm M}}_+}({{\rm d}}) 
\right\}\nonumber\\
{\mathbb P}_{-}({\rm M}_-)&=&\left\{({{\rm d}}, {{\rm f}})\in {\mathbb R}^2:\ {{\rm d}}\ge 0\,,\qquad   {{\rm f}}
\ge-{\rm M}_-{{\rm d}}\,,\quad 
{{\rm f}}
\ne {{\rm f}}^{\rm s}_{{{\rm M}}_-}({{\rm d}}) 
\right\}
\end{eqnarray}
The periods ${\tau}_+({\rm M}_+)$, ${\tau}_-({\rm M}_-)$ in~\eqref{TaTb} are given by
 \small{\begin{eqnarray}\label{tau-}
{\tau}_-({\rm M}_-)&=&2\sqrt{\frac{2\rm d}{v_0}}\int_{-1}^{\min\{1, \beta_-\}}\frac{d\beta}{ \sqrt{(1-\beta^2)\big({\rm d}^2\beta^2-4{\rm d}{{\rm M}}_-\beta+4{\rm f}-{\rm d}^2\big)}}\ {\rm if}\ ({{\rm d}}, {{\rm f}})\in
{\mathbb P}_{-}({\rm M}_-)
\\
\label{tau+}{\tau}_+({\rm M}_+)&=&2\sqrt{\frac{2\rm d}{v_0}}\int_{\max\{\alpha_-, 1\}}^{\alpha_+}\frac{d\alpha}{ \sqrt{(\alpha^2-1)\big(-{\rm d}^2\alpha^2+4{\rm d}{{\rm M}}_+\alpha-4{\rm f}+{\rm d}^2\big)}}\ {\rm if}\ ({{\rm d}}, {{\rm f}})\in
{\mathbb P}_{+}({\rm M}_+)
\end{eqnarray}}
\!\!\!where
\begin{eqnarray}\label{alphapm}\alpha_\pm:=\frac{2}{{{\rm d}}}\left( {{\rm M}}_+\pm \sqrt{{{\rm M}}_+^2+\frac{{{\rm d}}^2}{4}-{{\rm f}}}\right)\,,\ \beta_-({{\rm d}}, {{\rm f}})=\left\{\begin{array}{llll}\frac{2}{{{\rm d}}}\left({{\rm M}}_-- \sqrt{{{\rm M}}_-^2+\frac{{{\rm d}}^2}{4}-{{\rm f}}}\right)\ &{\rm if}\ {{\rm f}}\le {{\rm M}}_-^2+\frac{{{\rm d}}^2}{4}\\
1&{\rm otherwise}\,.
\end{array}
\right.
\end{eqnarray}
Furthermore, the following identity holds:
\begin{eqnarray}\label{tau+tau-}{\tau}_+({\rm M})=\tau_-({\rm M})\quad {\rm if}\ ({{\rm d}}, {{\rm f}})\in {\mathbb P}({\rm M}, {\rm M})\,.\end{eqnarray}\end{proposition}
\vskip.1in
\noindent

\proof  Let $\Theta_0=0$. It is immediate to see that  values of ${\rm J}_0\ge 0$ do not give rise to periodic motions, so we restrict to ${\rm J}_0< 0$, and denote as $\overline{\mathbb P}_\pm:={\mathbb H}_\pm\cap \{{\rm J}_0< 0\}$. Using the parameters $({{\rm d}}, {{\rm f}})$ in~\eqref{new parameters}, we show that
 \begin{eqnarray}\label{H-}\overline{\mathbb P}_+=\left\{({{\rm d}}, {{\rm f}})\in {\mathbb R}^2:\ {{\rm d}}\ge 0\,,\   {{\rm f}}
\le {{\rm f}}^+_{{{\rm M}}_+}({{\rm d}}) 
\right\}\,,\ 
\overline{\mathbb P}_-=\left\{({{\rm d}}, {{\rm f}})\in {\mathbb R}^2:\ {{\rm d}}\ge 0\,,\  {{\rm f}}
\ge -{\rm M}_-{\rm d}
\right\}\,.
\end{eqnarray}
Conditions in~\eqref{Hill}  are that the systems
\begin{eqnarray*}\left\{
\begin{array}{lll}
\displaystyle{{\rm d}^2}\alpha^2-4  {{\rm M}}_+{{\rm d}}\alpha+4 {{\rm f}}-{{\rm d}^2}\le 0\\\\
\displaystyle \alpha\ge1
\end{array}
\right.\qquad \left\{
\begin{array}{lll}
\displaystyle{{\rm d}^2}\beta^2-4  {{\rm M}}_-{{\rm d}}\beta+4 {{\rm f}}-{{\rm d}^2}\ge 0\\\\
-1\le\beta\le1
\end{array}
\right.\end{eqnarray*}
have non--empty solutions for $\alpha$, $\beta$, respectively.
Introducing the notations
\begin{eqnarray}\label{P1}
&&\overline{\mathbb P}_{{\rm M}}^{(1)}:=\Big\{({{\rm d}}, {{\rm f}})\in {{\mathbb R}}^2:\ {{\rm d}}\ge 0\,,\ -{{\rm M}}{{\rm d}}\le{{\rm f}}\le{{\rm M}}{{\rm d}}\Big\}\nonumber\\
&&\overline{\mathbb P}_{{\rm M}}^{(2)}:=\Big\{({{\rm d}}, {{\rm f}})\in {{\mathbb R}}^2:\  {{\rm d}}\ge2{{\rm M}}\,,\ {{\rm M}}{{\rm d}}\le{{\rm f}}\le{{\rm M}}^2+\frac{{{\rm d}}^2}{4}\Big\}\nonumber\\
&&\overline{\mathbb P}_{{\rm M}}^{(3)}:=\Big\{({{\rm d}}, {{\rm f}})\in {{\mathbb R}}^2:\ 0\le{{\rm d}}\le 2{{\rm M}}\,,\ {{\rm M}}{{\rm d}}\le{{\rm f}}\le{{\rm M}}^2+\frac{{{\rm d}}^2}{4}\Big\}\nonumber\\
&&\overline{\mathbb P}_{{\rm M}}^{(4)}:=\Big\{({{\rm d}}, {{\rm f}})\in {{\mathbb R}}^2:\   {{\rm d}}\ge 0\,,\ {{\rm f}}\ge
{{\rm M}}^2+\frac{{{\rm d}}^2}{4}
\Big\}\nonumber\\
&&\overline{\mathbb P}_{{\rm M}}^{(5)}:=\Big\{({{\rm d}}, {{\rm f}})\in {{\mathbb R}}^2:\   {{\rm d}}\ge 0\,,\ {{\rm f}}\le
-{{\rm M}}{\rm d}
\Big\}
\end{eqnarray}
   \begin{figure}
 \begin{center}
 \begin{tikzpicture}
\draw[densely dashed,  thick] (2,6) parabola bend (0,5)  (4,9);
\draw[ultra thick] (2,6) parabola bend (0,5)  (2,6);
\draw[densely dashed,  thick] plot coordinates {(0,4) (2,6)};
\draw[ultra thick] plot coordinates {(2,6) (5,9)};
\draw[ultra thick]  plot coordinates {(0,4) (5,0)};
\draw [->] (-2,4) -- (6,4);
\draw [->] (0,0) -- (0,10);
\node at (5.5,4.3) {d};
\node at (-0.3,9.5) {f};
\node at (4.0,9.5) {${\rm f=f^{\rm s}_{\rm M}({\rm d})}$};
\node at (5.5,8.5) {${\rm f=f^+_{\rm M}({\rm d})}$};
\node at (5.5,0.5) {${\rm f=f^-_{\rm M}({\rm d})}$};
\node at (3,4.8) {$\overline{\mathbb P}^{(1)}_{\rm M}$};
\node at (4.2,8.7) {$\overline{\mathbb P}^{(2)}_{\rm M}$};
\node at (0.4,4.8) {$\overline{\mathbb P}^{(3)}_{\rm M}$};
\node at (2.5,8.7) {$\overline{\mathbb P}^{(4)}_{\rm M}$};
\node at (0.4,3) {$\overline{\mathbb P}^{(5)}_{\rm M}$};
\end{tikzpicture}
  \end{center}
  \caption{Graphical representation of the lines ${\rm f}={\rm f}^-_{\rm M}({\rm d})$ (thick), ${\rm f}={\rm f}^+_{\rm M}({\rm d})$ (thick), ${\rm f}={\rm f}^{\rm s}_{\rm M}({\rm d})$ (dashed) and of the domains $\overline{\mathbb P}_{\rm M}^{(i)}$, in the plane ${\rm (d, f)}$.}\label{domains}
 \end{figure}
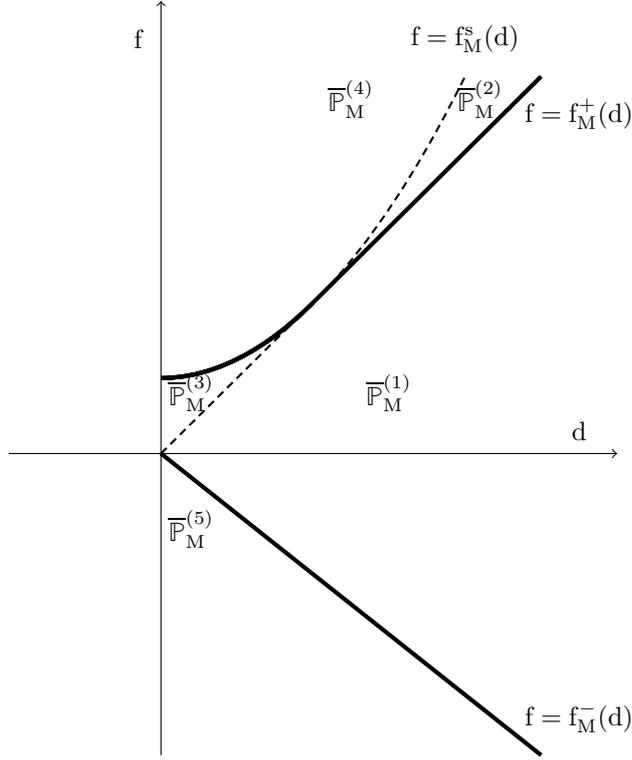
\!(see Fig.~\ref{domains}) one sees that
the solutions of such systems are
\begin{eqnarray}\label{solutions}
&&\left\{\begin{array}{llll}
\alpha_-\le \alpha\le \alpha_+\phantom{\ \&\ \beta_+\le \beta\le 1}\quad &{\rm if}\ ({\rm d}, {\rm f})\in \overline{\mathbb P}^{(3)}_{\rm M_+}\\
1\le \alpha\le\alpha_+\quad &{\rm if}\ ({\rm d}, {\rm f})\in  \overline{\mathbb P}^{(1)}_{\rm M_+}\cup \overline{\mathbb P}^{(5)}_{\rm M_+}\\
\emptyset&{\rm otherwise}\end{array}
\right.\nonumber\\
&&\left\{\begin{array}{llll}
-1\le \beta\le 1\phantom{\ \&\ \beta_+\le \beta\le 1}\quad &{\rm if}\ ({\rm d}, {\rm f})\in  \overline{\mathbb P}^{(3)}_{\rm M_-}\cup\overline{\mathbb P}^{(4)}_{\rm M_-}
\\
-1\le \beta\le \beta_-\phantom{\ \&\ \beta_+\le \beta\le 1}\quad &{\rm if}\ ({\rm d}, {\rm f})\in \overline{\mathbb P}^{(1)}_{\rm M_-}\\
-1\le\beta\le\beta_-\ \&\ \beta_+\le \beta\le 1\quad &{\rm if}\ ({\rm d}, {\rm f})\in \overline{\mathbb P}^{(2)}_{\rm M_-}\\
\emptyset&{\rm otherwise}
\end{array}
\right.\end{eqnarray}
where
\begin{eqnarray*}
\alpha_\pm:=\frac{2}{{{\rm d}}}\left( {{\rm M}}_+\pm \sqrt{{{\rm M}}_+^2+\frac{{{\rm d}}^2}{4}-{{\rm f}}}\right)\,,\ \beta_\pm=\frac{2}{{{\rm d}}}\left({{\rm M}}_-\pm \sqrt{{{\rm M}}_-^2+\frac{{{\rm d}}^2}{4}-{{\rm f}}}\right)\,.\end{eqnarray*}
Taking the union of the values of $(\rm d, f)$ which provide 
 non--empty solutions separately for $\alpha$ and $\beta$, we obtain the sets \begin{eqnarray*}\overline{\mathbb P}_+=\overline{\mathbb P}_{{\rm M}_+}^{(1)}\cup \overline{\mathbb P}_{{\rm M}_+}^{(3)}\cup \overline{\mathbb P}_{{\rm M}_+}^{(5)}\,,\qquad \overline{\mathbb P}_-=\overline{\mathbb P}_{{\rm M}_-}^{(1)}\cup \overline{\mathbb P}_{{\rm M}_-}^{(2)}\cup \overline{\mathbb P}_{{\rm M}_-}^{(3)}\cup \overline{\mathbb P}_{{\rm M}_-}^{(4)}\,.\end{eqnarray*} These are just the sets  in~\eqref{H-}.
 The level curves in~\eqref{eq for alpha},~\eqref{eq for beta} for $({\rm d}, {\rm f})\in \overline{\mathbb P}_+$, $\overline{\mathbb P}_-$, respectively, are
\begin{eqnarray}\label{new polynomials}&&\alpha'=\pm\sqrt{\frac{v_0}{2\rm d}}\sqrt{(\alpha^2-1)\big(-{\rm d}^2\alpha^2+4{\rm d}{{\rm M}}_+\alpha-4{\rm f}+{\rm d}^2\big)}\\
&&\label{new polynomials1}\beta'=\pm\sqrt{\frac{v_0}{2\rm d}} \sqrt{(1-\beta^2)\big({\rm d}^2\beta^2-4  {{\rm M}}_-{{\rm d}}\beta+4 {{\rm f}}-{\rm d}^2\big)}\,.
\end{eqnarray}
They  intersect the horizontal axes, in both $(\alpha, \alpha')$ and $(\beta, \beta')$ planes at  the end--points of the non--empty solutions for $\alpha$, $\beta$ (respectively) in~\eqref{solutions}.
The sets ${\mathbb P}_\pm$ is then the subset of $\overline{\mathbb P}_\pm$ such that the polynomials
under the square roots in ~\eqref{new polynomials},~\eqref{new polynomials1} have not double roots at such end--points.  
But such double roots occur if an only if  ${\rm f}= {{\rm f}}^{\rm s}_{{{\rm M}}_+}({{\rm d}})$ (for the $(\alpha, \alpha')$ level curve), or ${\rm f}= {{\rm f}}^{\rm s}_{{{\rm M}}_-}({{\rm d}})$ (for the $(\beta, \beta')$ level curve). Then ${\mathbb P}_\pm=\overline{\mathbb P}_\pm\setminus\{{\rm f}= {{\rm f}}^{\rm s}_{{{\rm M}}_\pm}({{\rm d}})\}$. These are precisely the sets in~\eqref{P0-}. 
At this point, the formulae in~\eqref{tau-},~\eqref{tau+} immediately follow, with the observation that
if $({{\rm d}}, {{\rm f}})\in \overline{\mathbb P}_{\rm M_-}^{(2)}\cap {\mathbb P}$, two periodic orbits for $\beta$
arise (one with $-1\le\beta\le \beta_-$, another with $\beta_+\le\beta\le 1$), but they have the same period, by an easy consequence of Cauchy Theorem.
Also the equality in~\eqref{tau+tau-} follows from the Cauchy Theorem. It is crucial to observe at this respect that, while $({{\rm d}}, {{\rm f}})\in {\mathbb P}({\rm M}, {\rm M})$, the roots of the polynomial ${\rm d}^2x^2-4{\rm d}{{\rm M}}x+4{\rm f}-{\rm d}^2$ are real--valued, and never fall in $(-\infty, -1)$. $\qquad \square$
 
 \vskip.1in
 \noindent
 The sets ${\mathbb P}_+$, ${\mathbb P}_+$ in~\eqref{P0-} are strictly larger than ${\mathbb P}({\rm M}_+, {\rm M}_-)$ in~\eqref{P0-OLD}. We show that if $({\rm d}, {\rm f})\in {\mathbb P}({\rm M}_+, {\rm M}_-)$ then the formulae for $\tau_-$, $\tau_+$ obtained through ~\eqref{tau-},~\eqref{tau+} coincide with the ones obtained through Equation~\eqref{TupFORMULA} and Proposition~\ref{thm: classical periods}.
  
\proof{\bf of Proposition~\ref{thm: classical periods}}
We first check the formula for $\tau_+$. At this purpose, observe that, as ${\rm M}_+\ge {\rm M}_-$, then  ${\mathbb P}({\rm M}_+, {\rm M}_-)\subset {\mathbb P}({\rm M}_+, {\rm M}_+)$. Then we can use~\eqref{tau+tau-} with ${\rm M}={\rm M}_+$ and~\eqref{tau-} with ${\rm M}_-$ replaced by ${\rm M}_+$. We obtain
\begin{eqnarray*}{\tau}_+({\rm M}_+)=\tau_-({\rm M}_+)=2\sqrt{\frac{2\rm d}{v_0}}\int_{-1}^{\min\{1, \beta^+_-\}}\frac{d\beta}{ \sqrt{(1-\beta^2)\big({\rm d}^2\beta^2-4{\rm d}{{\rm M}}_+\beta+4{\rm f}-{\rm d}^2\big)}}\end{eqnarray*}
for all $({{\rm d}}, {{\rm f}})\in {\mathbb P}({\rm M}_+, {\rm M}_-)$, with (observing that
the occurrence
${{\rm f}}> {{\rm M}}_+^2+\frac{{{\rm d}}^2}{4}$
never happens for $({{\rm d}}, {{\rm f}})\in {\mathbb P}({\rm M}_+, {\rm M}_-)$)
\begin{eqnarray*}\beta_-^+=\frac{2}{{{\rm d}}}\left({{\rm M}}_+- \sqrt{{{\rm M}}_+^2+\frac{{{\rm d}}^2}{4}-{{\rm f}}}\right)\end{eqnarray*}
being the minimal root of the polynomial ${\rm d}^2\beta^2-4{\rm d}{{\rm M}}_+\beta+4{\rm f}-{\rm d}^2$, hence real for all $({{\rm d}}, {{\rm f}})\in {\mathbb P}({\rm M}_+, {\rm M}_-)$.
It is immediate to check that $\beta_-^+> 1$ if and only if $({{\rm d}}, {{\rm f}})\in {\mathbb P}({\rm M}_+, {\rm M}_-)\cap {{\mathbb Q}}^\uparrow_{\rm M_+}$, whence we have nothing else to prove for $\tau_+$. 
Now we check the formula for ${\tau}_-$. The definition of $\beta_-$ in~\eqref{alphapm}
implies that the supremum in the integral~\eqref{tau-} is strictly less than $1$, and coincides with the minimum root of ${\rm d}^2x^2-4{\rm d}{{\rm M}}x+4{\rm f}-{\rm d}^2$ if and only if $({{\rm d}}, {{\rm f}})\in {\mathbb P}({\rm M}_+, {\rm M}_-)\cap {{\mathbb Q}}^\downarrow_{\rm M_-}\quad \square$.

\section{Dynamical part}\label{Proof of Theorem1.2}
The purpose of this and the next section is to provide a new formula 
for the function $\tau_{\rm M}$ in~\eqref{TupFORMULA}. Especially, we aim at simplifying its expression on ${{\mathbb Q}}^{\downarrow}_{\rm M}$.
\vskip.1in
\noindent
Define the function
\begin{eqnarray*}
{\rm t}_{\rm M}=\left\{\begin{array}{lll} \displaystyle{\rm t}^\uparrow_{\rm M}:=\sqrt{\frac{2\rm d}{v_0}}\int_{0}^{z_{\rm M}({\rm d}, {\rm f})}\frac{dz}{\sqrt{z(2-z)({\rm M}^2 z^2-2{\rm f} z+{\rm d}^2)}}\quad &{\rm if}\ ({{\rm d}}, {{\rm f}})\in
{{\mathbb Q}}^{\uparrow}_{\rm M}\\\\
\displaystyle{\rm t}^\downarrow_{\rm M}:=\sqrt{\frac{2\rm d}{v_0}}\int_{0}^{2}\frac{dz}{\sqrt{z(2-z)({\rm M}^2 z^2-2{\rm f} z+{\rm d}^2)}}
\quad &{\rm if}\ ({{\rm d}}, {{\rm f}})\in {{\mathbb Q}}^\downarrow_{\rm M}
\end{array}
\right.
\end{eqnarray*}
where  $z_{\rm M}({\rm d}, {\rm f})$ is the minimum roots of ${\rm M}^2 z^2-2{\rm f} z+{\rm d}^2$, turning out to be real while $({\rm d}, {\rm f})\in {{\mathbb Q}}^\uparrow_{\rm M}$. We shall prove that
\begin{proposition}\label{prop: new period}
The function ${\rm t}_{\rm M}$ is well--defined on ${{\mathbb Q}}_{\rm M}$
and verifies\begin{eqnarray*}
{\rm t}_{\rm M}=\tau_{\rm M}\quad \ \forall\ {	\rm M}\ge 0\,,\ \forall\ ({\rm d}, {\rm f})\in {{\mathbb Q}}_{\rm M}
\end{eqnarray*}
\end{proposition}
Proposition~\ref{prop: new period} will be proved in two steps, the former developed in this section, the latter in the next one. Similarly as in the proof of Proposition~\ref{classical Jacobi periods}, we consider the following  subsets of ${{\mathbb Q}}_{\rm M}$
\begin{eqnarray}\label{P1NEW}
&&{\mathbb P}_{{\rm M}}^{(1)}:=\Big\{({{\rm d}}, {{\rm f}})\in {{\mathbb R}}^2:\ {{\rm d}}\ge 0\,,\ -{{\rm M}}{{\rm d}}<{{\rm f}}<{{\rm M}}{{\rm d}}\Big\}\nonumber\\
&&{\mathbb P}_{{\rm M}}^{(2)}:=\Big\{({{\rm d}}, {{\rm f}})\in {{\mathbb R}}^2:\  {{\rm d}}\ge2{{\rm M}}\,,\ {{\rm M}}{{\rm d}}<{{\rm f}}<{{\rm M}}^2+\frac{{{\rm d}}^2}{4}\Big\}\nonumber\\
&&{\mathbb P}_{{\rm M}}^{(3)}:=\Big\{({{\rm d}}, {{\rm f}})\in {{\mathbb R}}^2:\ 0<{{\rm d}}< 2{{\rm M}}\,,\ {{\rm M}}{{\rm d}}<{{\rm f}}<{{\rm M}}^2+\frac{{{\rm d}}^2}{4}\Big\}\nonumber\\
&&{\mathbb P}_{{\rm M}}^{(4)}:=\Big\{({{\rm d}}, {{\rm f}})\in {{\mathbb R}}^2:\   {{\rm d}}\ge 0\,,\ {{\rm f}}>
{{\rm M}}^2+\frac{{{\rm d}}^2}{4}
\Big\}
\end{eqnarray}
which are the inner parts of the ones in~\eqref{P1}, hence are open. 
The union the ${\mathbb P}_{{\rm M}}^{(j)}$'s above coincides with ${{\mathbb Q}}_{\rm M}$, apart for the boundaries of the ${\mathbb P}_{{\rm M}}^{(j)}$ belonging to ${{\mathbb Q}}_{\rm M}$, which may be neglected for the purpose of proving Proposition~\ref{prop: new period}.
Note that, in the case $\rm M=0$,  ${\mathbb P}^{(1)}_{0}={\mathbb P}^{(3)}_{0}=\emptyset$. 
In this section, we prove that \begin{eqnarray}\label{first equality}\tau_{\rm M}=\sqrt{\frac{2\rm d}{v_0}} \int_{0}^{2}\frac{d z}{\sqrt{{z(2-z)({\rm M}}^2z^2-2{{\rm f}}z+{{\rm d}}^2)}}\ \forall\ {\rm M}>0\,,\ \forall\ ({\rm d}, {\rm f})\in {\mathbb P}^{(1)}_{\rm M}
\end{eqnarray}
In the next Section~\ref{A simpler formula}, we shall prove
\begin{eqnarray}\label{second equality}\tau_{\rm M}=\sqrt{\frac{2\rm d}{v_0}} \int_{0}^{\min\{2\,,{z_{\rm M}({\rm d}, {\rm f})}\}}\frac{d z}{\sqrt{{z(2-z)({\rm M}}^2z^2-2{{\rm f}}z+{{\rm d}}^2)}}\ \forall\ {\rm M}\ge 0\,,\  \forall\ ({\rm d}, {\rm f})\in {\mathbb P}^{(2)}_{\rm M}\cup {\mathbb P}^{(3)}_{\rm M}\cup {\mathbb P}^{(4)}_{\rm M}
\end{eqnarray}
(this equality includes the case $\rm M=0$, where one has to keep in mind  ${\mathbb P}^{(3)}_{0}=\emptyset$). 
Using \begin{eqnarray*}\overline{{\mathbb P}^{(1)}\cup {\mathbb P}^{(2)}}=\overline{{{\mathbb Q}}_{\rm M}^{\downarrow}}\,,\ \overline{{\mathbb P}^{(3)}\cup {\mathbb P}^{(4)}}=\overline{{{\mathbb Q}}_{\rm M}^{\uparrow}}\,,\ {z_{\rm M}({\rm d}, {\rm f})}\ge 2\,,\ \forall\ ({\rm d}, {\rm f})\in {\mathbb P}^{(2)}_{\rm M}
\end{eqnarray*}
The ``bar'' denotes set closure,
 we have that~\eqref{first equality} and~\eqref{second equality} imply Proposition~\ref{prop: new period}.
\subsection{Set up}
We fix the following terms.

\begin{itemize}
\item[\tiny\textbullet]
The coordinates \begin{eqnarray*}(A({\mathbf y}, {\mathbf x}), B({\mathbf y}, {\mathbf x}), \Theta({\mathbf y}, {\mathbf x}), \alpha({\mathbf x}), \beta({\mathbf x}), \vartheta({\mathbf x}))\end{eqnarray*} described in the introduction will be referred to as {\it co--focal coordinates associated to $({\mathbf y}, {\mathbf x})$}.

\item[\tiny\textbullet]Time $\tau$--solutions\footnote{As a rule, when the discussion  refers to the general problem, we keep the initial parameters $({\rm J}_0, {\rm F}_0, \Theta_0)$. The parameters $({\rm d}, {\rm f})$ in~\eqref{new parameters} will be used only in the planar case ($\Theta_0=0$).}
\begin{eqnarray}\label{standard solution}
\left\{\begin{array}{lll}\displaystyle\tau\to\big(A(\tau, {\rm J}_0, {\rm F}_0, \Theta_0, \alpha_0),  \alpha(\tau, {\rm J}_0, {\rm F}_0, \Theta_0, \alpha_0)\big)\\\\
\displaystyle \tau\to\big(B(\tau, {\rm J}_0, {\rm F}_0, \Theta_0, \beta_0),  \beta(\tau, {\rm J}_0, {\rm F}_0, \Theta_0,  \beta_0)\big)
\end{array}
\right.\end{eqnarray}
of the Hamiltonians ${\rm J}_+$, ${\rm J}_-$ in~\eqref{Jpm} satisfying~\eqref{AB} and, moreover,
\begin{eqnarray*}
\left\{\begin{array}{lll}\displaystyle \alpha(0, {\rm J}_0, {\rm F}_0, \Theta_0, \alpha_0)=\alpha_0\\\\
\displaystyle \beta(0, {\rm J}_0, {\rm F}_0, \Theta_0,  \beta_0)=\beta_0\end{array}
\right.\end{eqnarray*}
are called {\it Jacobi solutions determined by} $({\rm J}_0, {\rm F}_0, \Theta_0, \alpha_0, \beta_0)$.
  \item[\tiny\textbullet] 
  The periods  of the Jacobi solutions~\eqref{standard solution}, namely, the functions  $\tau_+({\rm J}_0, {\rm F}_0, \Theta_0)$, $\tau_-({\rm J}_0, {\rm F}_0, \Theta_0)$ in~\eqref{TaTb}, will be called {\it Jacobi periods}.

\item[\tiny\textbullet] Time $t$--solutions 
\begin{eqnarray}\label{ABabsolutions}
\left\{\begin{array}{lll}\displaystyle t\to\big(\overline A(t, {\rm J}_0, {\rm F}_0, \Theta_0, \alpha_0, \beta_0),  \overline\alpha(t, {\rm J}_0, {\rm F}_0, \Theta_0, \alpha_0, \beta_0)\big)\\\\
\displaystyle  t\to\big(\overline B(t, {\rm J}_0, {\rm F}_0, \Theta_0, \alpha_0, \beta_0),  \overline\beta(t, {\rm J}_0, {\rm F}_0, \Theta_0,  \alpha_0, \beta_0)\big)
\end{array}
\right.\end{eqnarray}
of the Hamiltonian ${\rm J}$ in~\eqref{Jnew} verifying the second equality in~\eqref{Jnew}, the equalities in~\eqref{AB} and, moreover,
\begin{eqnarray*}
\left\{\begin{array}{lll}\displaystyle \overline\alpha(0, {\rm J}_0, {\rm F}_0, \Theta_0, \alpha_0, \beta_0)=\alpha_0\\\\
\displaystyle \overline\beta(0, {\rm J}_0, {\rm F}_0, \Theta_0,  \alpha_0, \beta_0)=\beta_0\end{array}
\right.\end{eqnarray*}
will be called  {\it natural solutions determined by} $({\rm J}_0, {\rm F}_0, \Theta_0, \alpha_0, \beta_0)$.
 \item[\tiny\textbullet]  The periods   of the natural solutions~\eqref{ABabsolutions} will be called {\it natural periods}, and denoted as $\overline{\rm t}_+({\rm J}_0, {\rm F}_0, \Theta_0)$, $\overline{\rm t}_-({\rm J}_0, {\rm F}_0, \Theta_0)$.
\end{itemize}
\subsection{Jacobi periods}
In this section we construct the natural solutions~\eqref{ABabsolutions} starting from the solutions of the Hamiltonian ${\rm J}$ in~\eqref{oldHam}, regarding the parameters $({\rm J}_0, {\rm F}_0, \Theta_0)$ as the fixed values of their corresponding first integrals. 

\vskip.1in
\noindent
As we shall exploit the ``Keplerian limit'' (namely, the limit of~\eqref{oldHam} when $\rm m\to 0$), we preliminarily shift the initial   coordinates $({\mathbf y}, {\mathbf x})$ in~\eqref{oldHam} as
\begin{eqnarray}\label{shift}({\mathbf y}, {\mathbf x})\to ({\mathbf y}, {\mathbf x}-{\mathbf v}_0)\,.\end{eqnarray}
This  carries ${\rm J}$  to 
\begin{eqnarray}\label{newJ}{\rm J}({\mathbf y}, {\mathbf x})=\frac{\|{\mathbf y}\|^2}{2}-\frac{{{\rm M}}}{\|{\mathbf x}\|}-\frac{{\rm m}}{\|{\mathbf x}-2{\mathbf v}_0\|}\,.\end{eqnarray}
The function  $\Theta$ in~\eqref{m3}  is left unvaried by the transformation~\eqref{shift}. We let
\begin{eqnarray}\label{FFF}{\rm F}({\mathbf y}, {\mathbf x}, {{\rm M}}, {\rm m})=\|{\mathbf G}({\mathbf y}, {\mathbf x})\|^2-2{\mathbf v}_0\cdot\left(
{\mathbf y}\times{\mathbf G}({\mathbf y}, {\mathbf x})
-{{\rm M}}\frac{{\mathbf x}}{\|{\mathbf x}\|}+{\rm m}\frac{{\mathbf x}-2{\mathbf v}_0}{\|{\mathbf x}-2{\mathbf v}_0\|}\right)\end{eqnarray}
where
\begin{eqnarray}\label{G}{\mathbf G}({\mathbf y}, {\mathbf x}):={\mathbf x}\times {\mathbf y}\end{eqnarray}
 the angular\footnote{Note that, because of the shift~\eqref{shift}, ${\mathbf G}({\mathbf y}, {\mathbf x})$ in~\eqref{G} is related to ${\mathbf M}({\mathbf y}, {\mathbf x})$ in~\eqref{M} via
\begin{eqnarray*}{\mathbf G}={\mathbf M}+{\mathbf M}_0\quad {\rm with}\quad {\mathbf M}_0:={\mathbf v}_0\times {\mathbf y}\end{eqnarray*}} momentum. The function ${\rm F}({\mathbf y}, {\mathbf x}, {{\rm M}}, {\rm m})$ in~\eqref{FFF} will be referred to as
 {\it Euler Integral}. 
 
 \begin{proposition}\label{batxbary}
The  couple of functions~\eqref{ABabsolutions}
is a {\it natural  solution  determined by $({\rm J}_0, {\rm F}_0, \Theta_0, \alpha_0, \beta_0)$} if and only if there exists a solution \begin{eqnarray*} t\to \big(\overline{\mathbf y}(t, {\rm J}_0, {\rm F}_0, \Theta_0, \alpha_0, \beta_0)\,,\ \overline{\mathbf x}(t, {\rm J}_0, {\rm F}_0, \Theta_0, \alpha_0, \beta_0)\big)\end{eqnarray*}
of the Hamilton equations of ${\rm J}$ in~\eqref{newJ}, with initial datum $({\mathbf y}_0, {\mathbf x}_0)$ such that
\begin{eqnarray}\label{old system}\left\{
\begin{array}{lll}
\displaystyle{\rm J}({\mathbf y}_0, {\mathbf x}_0, {{\rm M}}, {\rm m})={\rm J}_0\\\\
\displaystyle{\rm F}({\mathbf y}_0, {\mathbf x}_0, {{\rm M}}, {\rm m})={\rm F}_0\\\\
\displaystyle\Theta({\mathbf y}_0, {\mathbf x}_0)=\Theta_0\\\\
\displaystyle \alpha({\mathbf x}_0)=\alpha_0\\\\
\beta({\mathbf x}_0)=\beta_0\,.
\end{array}
\right.
\end{eqnarray}
such that~\eqref{ABabsolutions} are the co--focal coordinates  associated to 
\begin{eqnarray*}\Big(\overline{\mathbf y}\big(t, {\rm J}_0, {\rm F}_0, \Theta_0, \alpha_0, \beta_0\big)\,,\  \overline{\mathbf x}\big(t, {\rm J}_0, {\rm F}_0, \Theta_0, \alpha_0, \beta\big)-{\mathbf v}_0\Big)\,.\end{eqnarray*}\end{proposition}
The proof of Proposition~\ref{batxbary} is standard, hence is deferred to Section~\ref{proof of Ham inter}.

\vskip.1in 
\noindent
In view of the definitions above, we rewrite the relation in~\eqref{dtaudt} between the times $t$ and $\tau$ as
\begin{eqnarray}\label{change of time}\tau=\int_0^t\frac{dt'}{\big\|\overline{\mathbf x}(t', {\rm J}_0, {\rm F}_0, \Theta_0, \alpha_0, \beta_0)\big\|\big\|\overline{\mathbf x}(t', {\rm J}_0, {\rm F}_0, \Theta_0, \alpha_0, \beta_0)-2{\mathbf v}_0\big\|}\end{eqnarray}
The  formula~\eqref{change of time} allows to establish the following relation between natural and Jacobi periods:
\begin{eqnarray}\label{Tab}
\tau_\pm({\rm J}_0, {\rm F}_0, \Theta_0)=\int_0^{\overline{\rm t}_\pm({\rm J}_0, {\rm F}_0, \Theta_0)}\ \frac{dt'}{\big\|\overline{\mathbf x}(t', {\rm J}_0, {\rm F}_0, \Theta_0, \alpha_0, \beta_0)\big\|\big\|\overline{\mathbf x}(t', {\rm J}_0, {\rm F}_0, \Theta_0, \alpha_0, \beta_0)-2{\mathbf v}_0\big\|}\end{eqnarray}
for any $\overline{\mathbf x}(t, {\rm J}_0, {\rm F}_0, \Theta_0, \alpha_0, \beta_0)$ as in Proposition~\ref{batxbary}.

\subsection{Jacobi Kepler period }\label{Jacobi Kepler periods}
In this section we assume ${\rm M}>0$. We consider
 the 
 Kepler's Hamiltonian
  \begin{eqnarray}\label{kepler}{\rm K}=\frac{\|{\mathbf y}\|^2}{2}-\frac{{\rm M}}{\|{\mathbf x}\|}\end{eqnarray}
which we regard as the limiting value of the Hamiltonian~\eqref{newJ} when 
  ${\rm m}=0$.
  We shall denote as
  \begin{eqnarray*}(\overline{\mathbf y}^{{\rm Kep}}(t, {\rm J}_0, {\rm F}_0, \Theta_0, \alpha_0, \beta_0), \overline{\mathbf x}^{{\rm Kep}}(t, {\rm J}_0, {\rm F}_0, \Theta_0, \alpha_0, \beta_0))\end{eqnarray*}
  any solution of ~\eqref{kepler} verifying 
\begin{eqnarray}\label{E***}\left\{
 \begin{array}{lll}
\displaystyle{\rm K}(\overline{\mathbf y}^{{\rm Kep}}, \overline{\mathbf x}^{{\rm Kep}})=\frac{\|\overline{\mathbf y}^{{\rm Kep}}\|^2}{2}-\frac{{\rm M}}{\|\overline{\mathbf x}^{{\rm Kep}}\|}={\rm J}_0\\\\
\displaystyle {\rm E}(\overline{\mathbf y}^{{\rm Kep}}, \overline{\mathbf x}^{{\rm Kep}})=\|{\mathbf G}(\overline{\mathbf y}^{{\rm Kep}}, \overline{\mathbf x}^{{\rm Kep}})\|^2-2{\mathbf v}_0\cdot\left(
\overline{\mathbf y}^{{\rm Kep}}\times{\mathbf G}(\overline{\mathbf y}^{{\rm Kep}}, \overline{\mathbf x}^{{\rm Kep}})
-{{\rm M}}\frac{\overline{\mathbf x}^{{\rm Kep}}}{\|\overline{\mathbf x}^{{\rm Kep}}\|}\right)={\rm F}_0\\\\
\Theta(\overline{\mathbf y}^{{\rm Kep}}, \overline{\mathbf x}^{{\rm Kep}})={\mathbf G}(\overline{\mathbf y}^{{\rm Kep}}, \overline{\mathbf x}^{{\rm Kep}})\cdot {\mathbf i}=\Theta_0
 \end{array}
 \right.\end{eqnarray} As well--known, the Keplerian orbit $t\to\overline{\mathbf x}^{{\rm Kep}}(t, {\rm J}_0, {\rm F}_0, \Theta_0, \alpha_0, \beta_0)$ describes an ellipse, the semi--major axis, eccentricity, anomaly of perihelion of which will be denoted as $a$, $e$, $\omega$, respectively. The  periods of the Keplerian orbits, here denoted as ${\rm t}_{{\rm M}}^{{\rm Kep}}$, are well--known to depend only  on the semi--major axis $a$, or, equivalently,  on the value ${\rm J}_0$ of the energy:
   \begin{eqnarray*}{\rm t}_{{\rm M}}^{{\rm Kep}}({\rm J}_0)=2\pi \sqrt{\frac{a^3}{{\rm M}}}=2\pi {{\rm M}}\sqrt{\frac{1}{8|{\rm J}_0|^3 }}\end{eqnarray*}
   as $a=\frac{{\rm M}}{2|{\rm J}_0|}$.
With the terminology of the previous section, we have that 
 the natural periods $\overline {\rm t}_+$, $\overline {\rm t}_-$ for the Keplerian case coincide, and their common value is precisely ${\rm t}_{{\rm M}}^{{\rm Kep}}({\rm J}_0)$.  But then, also the Jacobi periods  ({\it Jacobi Kepler periods}, hereafter)  coincide. They will be denoted as $\tau^{{\rm Kep}}_{{\rm M}}$ and, through~\eqref{Tab}, are given by
 \begin{eqnarray}\label{TabNEW}
\tau^{{\rm Kep}}_{{\rm M}}({\rm J}_0, {\rm F}_0, \Theta_0)=\int_0^{{\rm t}^{{\rm Kep}}_{{\rm M}}({\rm J}_0)}\ \frac{dt'}{\big\|\overline{\mathbf x}^{{\rm Kep}}(t', {\rm J}_0, {\rm F}_0, \Theta_0, \alpha_0, \beta_0)\big\|\big\|\overline{\mathbf x}^{{\rm Kep}}(t', {\rm J}_0, {\rm F}_0, \Theta_0, \alpha_0, \beta_0)-2{\mathbf v}_0\big\|}\end{eqnarray}
It is convenient to
use the eccentric anomaly $\xi$ in the integral~\eqref{TabNEW}, as it is related to $t'$ via:
\begin{eqnarray*}\label{txi} \frac{dt'}{d\xi}=\sqrt{\frac{a}{{\rm M}}}\|\overline {\mathbf x}^{{\rm Kep}}\|
\end{eqnarray*}
This way, $\|\overline{\mathbf x}^{{\rm Kep}}\|$ is cancelled, and~\eqref{TabNEW} becomes

 \begin{eqnarray*}\label{general t} \tau^{{\rm Kep}}_{{\rm M}}({\rm J}_0, {\rm F}_0, \Theta_0)=\sqrt{\frac{a}{{\rm M}}}\int_0^{2\pi}\ \frac{d\xi}{\big\|\overline{\mathbf x}^{{\rm Kep}}(\xi, {\rm J}_0, {\rm F}_0, \Theta_0, \alpha_0, \beta_0)-2{\mathbf v}_0\big\|}
\end{eqnarray*}

\vskip.1in
\noindent
We now focus on the case $\Theta_0=0$.  Using the orbital elements $a$, $e$, $\omega$, we  obtain
\begin{eqnarray}\label{tKepOLD}
\tau^{{\rm Kep}}_{{\rm M}}({\rm J}_0, {\rm F}_0)=\sqrt{\frac{a}{{\rm M}}}\int_0^{2\pi}\frac{d\xi}{\sqrt{a^2(1-e\cos\xi)^2-4 a v_0 \big((\cos\xi-e)\cos\omega-\sqrt{1-e^2}\sin\xi\sin \omega\big)+4v_0^2
}}
\end{eqnarray}
where $a$, $e$ and $\omega$ are to be regarded as functions of $({\rm J}_0, {\rm F}_0)$. To this end, we rewrite two first equations in~\eqref{E***}  using $a$, $e$ and $\omega$:
\begin{eqnarray*}
\left\{\begin{array}{lll}
\displaystyle-\frac{{\rm M}}{2a}={\rm J}_0\\\\
\displaystyle{{\rm M}}\left(a (1-e^2)-2v_0 e \cos\omega\right)={\rm F}_0
\end{array}
\right.
\end{eqnarray*}
Switching, via~\eqref{new parameters}, to $({\rm d}, {\rm f})$,
they become
\begin{eqnarray}\label{substitutions}
\left\{\begin{array}{lll}
\displaystyle a=\frac{2v_0\rm M}{\rm d}\\\\
\displaystyle e^2+\frac{{{\rm d}}}{{\rm M}} e \cos\omega=1-\frac{{{\rm f}}}{{{\rm M}}^2}
\end{array}
\right.
\end{eqnarray}
Observe that we have only two parameters, $({\rm d}, {\rm f})$, for  the three quantities $a$, $e$ and $\omega$. This will enable us to fix one among $(e, \omega)$ arbitrarily. 
At this respect, the second equation above has non--empty solutions for $(e, \omega)\in [0, 1)\times{{\mathbb T}}$ for any choice of $({{\rm d}}, {{\rm f}})$ in the set
\begin{eqnarray}\label{P0-KEPOLD} \overline{\mathbb P}^{\rm Kep}_{{\rm M}}:=\left\{({{\rm d}}, {{\rm f}}):\ {{\rm d}}\ge0\,,\qquad {{\rm f}}^-_{{\rm M}}({{\rm d}})\le {{\rm f}}
\le {{\rm f}}^+_{{\rm M}}({{\rm d}})\right\}\end{eqnarray}
with ${{\rm f}}^-_{{\rm M}}({{\rm d}})$, ${{\rm f}}^+_{{\rm M}}({{\rm d}})$ the functions we have already encountered in~\eqref{FmaxFsing}.
With
\begin{eqnarray}\label{P0-KEP} 
{\mathbb P}^{\rm Kep}_{{\rm M}}:=\overline{\mathbb P}^{\rm Kep}_{{\rm M}}\setminus\{ {{\rm f}}= {{\rm f}}^{\rm s}_{{\rm M}}({{\rm d}})\}
\end{eqnarray}
we have
\begin{proposition}\label{Keplerian periods}
$\tau^{\rm Kep}_{{\rm M}}({{\rm d}}, {{\rm f}})$ is bounded if and only if $({{\rm d}}, {{\rm f}})\in {\mathbb P}^{\rm Kep}_{{\rm M}}$.
Moreover, it is given by
\begin{eqnarray}\label{T+-}
\footnotesize
&&\tau^{\rm Kep}_{{\rm M}}({{\rm d}}, {{\rm f}})\nonumber\\
&&={\sqrt{\frac{\rm d}{2v_0}}}\int_0^{2\pi}\frac{d\xi}{\sqrt{{\rm M}^2(1-e\cos\xi)^2-2 {{\rm M}}{{\rm d}}\big((\cos\xi-e)\cos\omega-\sqrt{1-e^2}\sin\xi\sin \omega\big)+{{\rm d}}^2
}}\nonumber\\
\end{eqnarray}
for all $({{\rm d}}, {{\rm f}})\in {\mathbb P}^{\rm Kep}_{{\rm M}}$. Here, $(e, \omega)$ is any couple in $ [0, 1)\times {{\mathbb T}}$ solving \begin{eqnarray}\label{equation for e omega}
e^2+\frac{{{\rm d}}}{{\rm M}} e \cos\omega=1-\frac{{{\rm f}}}{{{\rm M}}^2}
\end{eqnarray} 
\end{proposition}
\proof 
The formula in~\eqref{T+-} is the same as in~\eqref{tKepOLD}, after the substituting $a$  via~\eqref{substitutions}.
For any $({{\rm d}}, {{\rm f}})\in \overline{\mathbb P}^{\rm Kep}_{{\rm M}}$, the couple
\begin{eqnarray*}(e_*, \omega_*)=\left(\frac{\widehat e}{{\rm M}}\,,\  {\cos^{-1}\sigma}\right)\end{eqnarray*} 
with \begin{eqnarray*}\widehat e:=\sigma\left(-\frac{{{\rm d}}}{2}+\sqrt{{{\rm M}}^2+\frac{{{\rm d}}^2}{4}-{{\rm f}}}\right)\qquad \sigma:=
\left\{
\begin{array}{lll}
\displaystyle+1\quad &{\rm if}\  -{{\rm M}}{{\rm d}}\le {{\rm f}}<{\rm M}^2\\\\
-1\quad &{\rm if}\   {\rm M}^2\le {{\rm f}}<{\rm f}^+_{{\rm M}}({{\rm d}})
\end{array}
\right.\end{eqnarray*}
verifies $(e_*, \omega_*)\in [0, 1)\times {\mathbb T}$ and solves~\eqref{equation for e omega}. Inserting such values in~\eqref{T+-}, the term with $\sin\xi$ is cancelled, and one obtains the formula
 \begin{eqnarray}\label{Tuparrow}
\tau^{\rm Kep}_{{\rm M}}({{\rm d}}, {{\rm f}})=\sqrt{\frac{\rm d}{2v_0}} \int_0^{2\pi}\frac{d \xi}{\sqrt{\widehat e^2\cos^2\xi^2-2{\rm M}(\widehat e+\sigma{{\rm d}})\cos\xi+{\rm M}^2+2\sigma\widehat e{{\rm d}}+{{\rm d}}^2}}
\end{eqnarray}
Now, the integral in~\eqref{Tuparrow} is (real--valued and) finite provided that the polynomial
\begin{eqnarray*} P(z):=\widehat e^2z^2-2{\rm M}(\widehat e+\sigma{{\rm d}})z+{\rm M}^2+2\sigma\widehat e{{\rm d}}+{{\rm d}}^2
\end{eqnarray*}
has no roots in ${\pm 1}$ and no double roots in $(-1, 1)$. This  defines the domain~\eqref{P0-KEP}. $\quad \square$

\vskip.1in
\noindent
The domains ${\mathbb P}^{\rm Kep}_{{\rm M}_-}$, ${\mathbb P}^{\rm Kep}_{{\rm M}_+}$ defined via~\eqref{P0-KEPOLD},~\eqref{P0-KEP}  are strictly smaller\footnote{Observe however that ${\mathbb P}^{\rm Kep}_{{\rm M}_+}$ is larger than ${\mathbb P}({\rm M}_+, {\rm M}_-)$ in~\eqref{P0-OLD}. This already allows us to identify $\tau_+=\tau^{\rm Kep}_{\rm M_+}$ for all $({\rm d}, {\rm f})\in{\mathbb P}({\rm M}_+, {\rm M}_-)$, even though a similar formula for $\tau_-$ is still missing at this stage.} than the domains ${\mathbb P}_-$, ${\mathbb P}_+$ in~\eqref{P0-}.
By construction, the periods  $\tau_-$, $\tau_{+}$ in~\eqref{tau-},~\eqref{tau+} coincide with $\tau^{\rm Kep}_{\rm M_+}$, $\tau^{\rm Kep}_{\rm M_-}$ limited to the values of $({\rm d}, {\rm f})\in {\mathbb P}^{\rm Kep}_{{\rm M}_+}$, ${\mathbb P}^{\rm Kep}_{{\rm M}_-}$, respectively.
Notice however that the formula~\eqref{Tuparrow} exhibited inside the proof of Proposition~\ref{Keplerian periods} is involved, and seems to be useless in order to prove Theorem~\ref{extendperiods}. However, restricting $({\rm d}, {\rm f})$ to the domains ${\mathbb P}_{\rm M_-}^{(1)}$, ${\mathbb P}_{\rm M_+}^{(1)}$ defined in~\eqref{P1NEW}, 
and moreover using -- as in  the proof of Proposition~\ref{Keplerian periods} -- the remarked freedom of choice of one among $(e, \omega)$ in equation~\eqref{equation for e omega}, 
we obtain a simpler formula, as we show in the next section.

\subsection{Proof of~\eqref{first equality}}
When $({{\rm d}}, {{\rm f}})$ takes values in the set ${\mathbb P}^{(1)}_{{\rm M}}$ in~\eqref{P1},
 the limiting couple
\begin{eqnarray*}(e^*, \omega^*)=\left(1, -\cos^{-1}\frac{{{\rm f}}}{{{\rm M}}{{\rm d}}}\right)\end{eqnarray*}
solves Equation~\eqref{equation for e omega}, and provides (through Proposition~\ref{Keplerian periods}) the expression \begin{eqnarray*}
{\rm t}^{*}_{{\rm M}}({{\rm d}}, {{\rm f}})&=&\sqrt{\frac{\rm d}{2v_0}} \int_{0}^{2\pi }\frac{d \xi}{\sqrt{{{\rm M}}^2(1-\cos\xi)^2-2{{\rm f}}(1-\cos\xi)+{{\rm d}}^2}}\nonumber\\
&=&\sqrt{\frac{2\rm d}{v_0}}\int_{0}^{2}\frac{d z}{\sqrt{{z(2-z)({\rm M}}^2z^2-2{{\rm f}}z+{{\rm d}}^2)}}\end{eqnarray*}
having let $z:=1-\cos\xi$ in the half--period.
This proves ~\eqref{first equality}.

\subsection{Proof of Proposition~\ref{batxbary}}\label{proof of Ham inter}
{The proof of the following lemma is omitted:}
\begin{lemma}
Using the coordinates $(A, B, \Theta, \alpha, \beta, \vartheta)$ described in Section~\ref{Purpose of the paper}, the {functions ${\rm J}$ and ${\rm F}$ in~\eqref{newJ},~\eqref{FFF} have the expressions} {\begin{eqnarray}\label{Enew1}
\left\{\begin{array}{llll}
{\rm J}(A, B, \Theta, \alpha, \beta, {\rm M}, {\rm m})=\frac{1}{v_0^2(\alpha^2-\beta^2)}\big({\rm J}_-(B, \beta, 0, \Theta, {\rm M}_-)+{\rm J}_+(A, \alpha, 0, \Theta, {\rm M}_+)\big)\\\\
{\rm F}(A, B, \Theta, \alpha, \beta, {\rm M}, {\rm m})=2\frac{\alpha^2-1}{\alpha^2-\beta^2}{\rm J}_-(B, \beta, 0, \Theta, {\rm M}_-)-2\frac{1-\beta^2}{\alpha^2-\beta^2}{\rm J}_+(A, \alpha, 0, \Theta, {\rm M}_+)
\end{array}
\right.
\end{eqnarray}}
where  
{${\rm J}_+(A, \alpha, {{\rm J}_0, \Theta} , {{\rm M}})$, ${\rm J}_-(B, \beta, \Theta, {\rm J}_0 , {{\rm M}})$}
are as in~\eqref{Jpm}, {${\rm M}_+$, ${\rm M}_-$ as in~\eqref{M+-}}\,.
\end{lemma}

\proof {\bf of Proposition~\ref{batxbary}} We look at the system~\eqref{old system} using the co--focal coordinates. Namely, we look at motions of ${\rm J}(A, B, \Theta, \alpha, \beta, {{\rm M}}, {\rm m})$ in~\eqref{Jnew} with an initial datum $(A(0), B(0), \Theta(0), \alpha(0), \beta(0))$ such that
\begin{eqnarray}\label{level surface cofocal}\left\{
\begin{array}{lll}
\displaystyle{\rm J}(A(0), B(0), \Theta(0), \alpha(0), \beta(0), {{\rm M}}, {\rm m})={\rm J}_0\\\\
\displaystyle{\rm F}(A(0), B(0), \Theta(0), \alpha(0), \beta(0), {{\rm M}}, {\rm m})={\rm F}_0\\\\
\displaystyle\Theta(0)=\Theta_0\\\\
\displaystyle \alpha(0)=\alpha_0\\\\
\beta(0)=\beta_0\,.
\end{array}
\right.
\end{eqnarray} 
{Using~\eqref{Enew1}, we see that, if $(\alpha, \beta)\ne (1, \pm 1)$, the system~\eqref{level surface cofocal} is equivalent to
\begin{eqnarray*}
\left\{\begin{array}{llll}
\displaystyle{\rm J}_+(A, \alpha, {0, \Theta} , {{\rm M}_+})=v_0^2(\alpha^2-1){\rm J}_0-\frac{{\rm F}_0}{2}\\\\
\displaystyle{\rm J}_-(B, \beta, 0, \Theta , {{\rm M}_-})
=v_0^2(1-\beta^2){\rm J}_0+\frac{{\rm F}_0}{2}
\end{array}
\right.
\end{eqnarray*}
This is the same as~\eqref{AB}, as ${\rm J}_+(A, \alpha, {0, \Theta} , {{\rm M}_+})-v_0^2(\alpha^2-1){\rm J}_0={\rm J}_+(A, \alpha, {{\rm J}_0, \Theta} , {{\rm M}_+})$, and ${\rm J}_-(B, \beta, {0, \Theta} , {{\rm M}_-})-v_0^2(1-\beta^2){\rm J}_0={\rm J}_-(B, \beta, {{\rm J}_0, \Theta} , {{\rm M}_-})$. $\quad \square$
}

\section{Analytic part}\label{A simpler formula}

In this section, we prove a result on elliptic integrals. For notices on elliptic integrals, we refer to~\cite{bianchi1898}.
\subsection{Equivalences of complex integrals}\begin{lemma}
Given $\varrho\in {\mathbb C}\setminus\{\pm 1\}$, the unique circle in ${\mathbb C}$ through $z=1$, $z=-1$ and $z=\varrho$ has equation
\begin{eqnarray}\label{Crho}
|z|^2+\left(\frac{1-|\varrho|^2}{\Im \varrho}\right)\Im z=1
\end{eqnarray}
\end{lemma}
\begin{definition}\rm We denote as:
\item[\tiny\textbullet] ${\mathcal C}(\varrho)$ the circle  in~\eqref{Crho}; 
\item[\tiny\textbullet] ${\mathcal S}(\varrho)$ the segment of  ${\mathcal C}(\varrho)$ where $\varrho$ does not belong to;
\item[\tiny\textbullet] ${\mathcal S}^*(\varrho):={\mathcal C}(\varrho)\setminus{\mathcal S}(\varrho)$;
\item[\tiny\textbullet] ${\mathcal A}(\varrho)$ the {closed} bounded region delimited by ${\mathcal S}(\varrho)$ and the interval $[-1,1]$ in the real axis;
\item[\tiny\textbullet] ${\mathcal A}^*(\varrho)$ the {closed} bounded region delimited by ${\mathcal S}^*(\varrho)$ and the interval $[-1,1]$ in the real axis.
\end{definition}

\begin{lemma}\label{Gammazeta}
Let $ \varrho\in {\mathbb C}\setminus\{\pm1\}$. The invertible transformation
 \begin{eqnarray*}
\Gamma_\varrho:\qquad z\in {\mathbb C}\cup \{\infty\}\to \tau=  \frac{1-z \varrho}{z-\varrho}\in {\mathbb C}\cup \{\infty\}\end{eqnarray*}
  sends $1\to 1$, $(-1)\to(-1)$, $\infty\to \varrho$ and, moreover, the interval $[-1,1]$ in the $z$--plane to the path $s_\varrho$ in the $\tau$--plane, given by:
  \begin{itemize}
\item[\rm (i)] if $\varrho\in {\mathbb R}$ and $|\varrho|>1:$ the interval $[-1,1]^+$ in the $z$ to itself
in the $\tau$--plane; the half--plane $\{\Im z>0\}$, $(\{\Im z<0\})$ to itself in the $\tau$--plane;
\item[\rm (ii)]  if $\varrho\in {\mathbb R}$, $|\varrho|<1:$  the interval $[-1,1]^+$ in the $z$--plane to
 $(-\infty, -1]^-\cup [1, +\infty)^-$ in the $\tau$--plane; the half--plane $\{\Im z>0\}$, $(\{\Im z<0\})$ to $\{\Im \tau<0\}$, $(\{\Im \tau>0\})$;
\item[\rm (iii)] the circular segment  ${\mathcal S}(-\varrho)$, run from $\tau=-1$ to $\tau=1$, if $\varrho\in {\mathbb C}\setminus{\mathbb R}$.\end{itemize}
\end{lemma}
\proof We prove only (iii), as (i) and (ii) are trivial. Let $t\in [-1,1]$. We prove that  $\tau=\Gamma_\varrho(t)\in {\mathcal C}(-\varrho)$. Writing $\varrho=a+\ii b$ with $a$, $b\in {\mathbb R}$, we have
\begin{eqnarray*}
|\tau|^2-\frac{1-|\varrho|^2}{\Im (\varrho)}\Im \tau&=&\left|\frac{1-t (a+\ii b)}{t-(a+\ii b)}\right|^2-\frac{1-a^2-b^2}{b}\Im \left(\frac{1-t (a+\ii b)}{t-(a+\ii b)}\right)\nonumber\\
&=&\frac{(1-ta)^2+t^2b^2}{(t-a)^2+b^2}-\frac{1-a^2-b^2}{b}\frac{(1-ta)b-tb(t-a)}{(t-a)^2+b^2}\nonumber\\
&=&\frac{(1-ta)^2+t^2b^2}{(t-a)^2+b^2}-\frac{(1-a^2-b^2)(1-t^2)}{(t-a)^2+b^2}\nonumber\\
&=&1\,.
\end{eqnarray*}
On the other hand, the image of $[-1,1]$ through $\Gamma_\varrho$ must be the circular segment of ${\mathcal C}(-\varrho)$ with extremes $\tau=-1$, $\tau=1$ and lying on the opposite side with respect to  $-\varrho$, as $-\varrho$ is the image of $\infty$. This is precisely ${\mathcal S}(-\varrho)$. $\quad\square$
\vskip.1in
\noindent
The main results of this section is the following one.  Notice that we do not assume that $P$  and $R$ have real--valued coefficients.
\begin{proposition}\label{Legendre polynomials2COR}
Let \begin{eqnarray*}P(z)=az^2-2bz+c\,,\qquad R(z)=\alpha z^2-2\beta z+\gamma\end{eqnarray*}
where $a$, $b$, $c$, $\alpha$, $\beta$, $\gamma\in {\mathbb C}$
 verify
\begin{eqnarray}\label{conditions2}
&&a\pm 2b+c\ne 0
\\
\label{conditions5}
&&{\rm if} \ a\ne 0\ {\rm and}\ b^2-ac=0\ {\rm then}\ \sqrt{\frac{c}{a}}\notin (-1, 1)\\
\label{conditions3}&&a-c=\alpha-\gamma\\ 
\label{conditions4}&&b^2-ac=\beta^2-\alpha\gamma\,.
\end{eqnarray}
Choose  one root $\varrho\in{\mathbb C}$ of the polynomial
\begin{eqnarray}\label{eqforrho}
Q(z):=(a-\alpha)z^2-2b z+ c+\alpha
\,.\end{eqnarray}
Then $\varrho\notin\{\pm 1\}$ and 
\begin{itemize}
\item[$a)$] the equality
\begin{eqnarray}\label{quadratic integrals}\int_{[-1, 1]}\frac{dz}{\sqrt{(1-z^2)P(z)}}&=&\int_{[-1, 1]}\frac{dz}{\sqrt{(1-z^2)R(z)}}\end{eqnarray}
hold if
 \begin{itemize}
\item[$a_1)$] $\varrho\in {\mathbb R}$, with $|\varrho|>1$;
\item[$a_2)$] $\varrho\in {\mathbb R}$, with $|\varrho|<1$ and at least one of the half--planes $\{\Im z>0\}$, $\{\Im z<0\}$ is free from
 roots of $P$ (of $R$);
\item[$a_3)$] $\varrho\in {\mathbb C}\setminus {\mathbb R}$ and no root of $P$ lies in ${\mathcal A}(\varrho)$ (no root of $R$ lies in ${\mathcal A}(-\varrho)$);
\end{itemize}
\item[$b)$] the equality
\begin{eqnarray*} \int_{[-1, 1]}\frac{dz}{\sqrt{(1-z^2)P(z)}}&=&-\int_{[-\infty, -1]\cup [1, +\infty)}\frac{dz}{\sqrt{(1-z^2)R(z)}}\end{eqnarray*}
hold if 
 \begin{itemize}
\item[$b_1)$]$\varrho\in {\mathbb R}$, with $|\varrho|>1$ and at least one of the half--planes $\{\Im z>0\}$, $\{\Im z<0\}$ is free from
 roots of $P$ (of $R$);
\item[$b_2)$]$\varrho\in {\mathbb R}$, with $|\varrho|<1$;
 \item[$b_3)$] $\varrho\in {\mathbb C}\setminus {\mathbb R}$ and no root of $P$ lies in ${\mathcal A}^*(\varrho)$ (equivalently, no root of $R$ lies in ${\mathcal A}^*(-\varrho)$).
 \end{itemize}
 \item[$c)$] the equality
\begin{eqnarray*} \int_{[-\infty, -1]\cup [1, +\infty)}\frac{dz}{\sqrt{(1-z^2)P(z)}}&=&\int_{[-\infty, -1]\cup [1, +\infty)}\frac{dz}{\sqrt{(1-z^2)R(z)}}
\end{eqnarray*}
hold if
 \begin{itemize}
\item[$c_1)$] $\varrho\in {\mathbb R}$, with $|\varrho|>1$;
\item[$c_2)$] $\varrho\in {\mathbb R}$, with $|\varrho|<1$ and at least one of the half--planes $\{\Im z>0\}$, $\{\Im z<0\}$ is free from
 roots of $P$ (of $R$);
\item[$c_3)$] $\varrho\in {\mathbb C}\setminus {\mathbb R}$ and no root of $P$ lies in ${\mathcal A}^*(1/\varrho)$ (no root of $R$ lies in ${\mathcal A}^*(-1/\varrho)$);
\end{itemize}
\item[$d)$] the equality
\begin{eqnarray*} \int_{[-\infty, -1]\cup [1, +\infty)}\frac{dz}{\sqrt{(1-z^2)P(z)}}&=&-\int_{[-1, 1]}\frac{dz}{\sqrt{(1-z^2)R(z)}}\end{eqnarray*}
holds if 
 \begin{itemize}
\item[$d_1)$]$\varrho\in {\mathbb R}$, with $|\varrho|>1$ and at least one of the half--planes $\{\Im z>0\}$, $\{\Im z<0\}$ is free from
 roots of $P$ (of $R$);
\item[$d_2)$]$\varrho\in {\mathbb R}$, with $|\varrho|<1$;
 \item[$d_3)$] $\varrho\in {\mathbb C}\setminus {\mathbb R}$ and no root of $P$ lies in ${\mathcal A}(1/\varrho)$ (equivalently, no root of $R$ lies in ${\mathcal A}(-1/\varrho)$).
 \end{itemize}
 \item[$e)$] If $\varrho\in {\mathbb R}$, the equality
\begin{eqnarray*} \int_{[x_0, 1]}\frac{dz}{\sqrt{(1-z^2)P(z)}}=\int_{[y_0(\varrho), 1]}\frac{dz}{\sqrt{(1-z^2)R(z)}}\end{eqnarray*}
also holds, for any {$x_0\in {\mathbb R}$}, and with $y_0(\varrho)=\Gamma_\varrho(x_0)$.
 \end{itemize}
\end{proposition}

\begin{remark}\rm
The reason why the sets ${\mathcal A}(\varrho)$, ${\mathcal A}^*(\varrho)$, ${\mathcal A}^*(1/\varrho)$, ${\mathcal A}(1/\varrho)$ appearing in assumptions $a_3)$, $b_3)$, $c_3)$, $d_3)$ are chosen to be closed is to avoid annoying questions of definition of the complex square root (as the proof uses arguments of holomorphic functions).
\end{remark}
 	
\proof {\bf of Proposition~\ref{Legendre polynomials2COR}}
We prove only the statements in $a)$,  as the proof of the statements in $b)$  is similar, 
$b)$, $c)$
 are  obtained from $a)$, $b)$ respectively, 
 changing coordinate $z=\zeta^{-1}$. Finally, statement in $e)$ a generalization of $a_1)$, $a_2)$.\\
 Let $\varrho\in{{\mathbb C}}$ solve~\eqref{eqforrho}. Due to condition~\eqref{conditions2}, we have $\varrho\in{\mathbb C}\setminus\{\pm 1\}$.
We check that
$\alpha$, $\beta$ and $\gamma$ verify
 \begin{eqnarray}\label{alphabetagamma}\alpha=\frac{a\varrho^2-2b \varrho+ c}{\varrho^2-1}\,,\quad \beta=\frac{b\varrho^2{-}(a+c) \varrho+b}{\varrho^2-1}\,,\quad \gamma=\frac{c\varrho^2-2b \varrho+a}{\varrho^2-1}\,.\end{eqnarray}
 Indeed, the first equality in~\eqref{alphabetagamma} is equivalent to~\eqref{eqforrho} combined with $\varrho\in{\mathbb C}\setminus\{\pm 1\}$. The second and the third equality 
 follow from it and ~\eqref{conditions3} and~\eqref{conditions4}.
Let now ${\mathcal I}_1$ denote the integral
at left hand side in~\eqref{quadratic integrals}.
We change variable, letting \begin{eqnarray*}z=\Gamma_{-\varrho}(\tau)=\frac{1+\tau\varrho}{\tau+\varrho}\end{eqnarray*}
where the coordinate $\tau$ runs the path  $s_\varrho:=\Gamma_\varrho([-1,1])$, as in the thesis of Lemma~\ref{Gammazeta}. The integral becomes
 \begin{eqnarray*}
 {\mathcal I}_1=\int_{s_\varrho }\frac{w(\tau)d\tau}{\sqrt{(1-\tau^2)
 \big(\alpha\tau^2-2\beta \tau+\gamma\big)
 w(\tau)^2}}=\int_{s_\varrho }\frac{d\tau}{\sqrt{(1-\tau^2)
 \big(\alpha\tau^2-2\beta \tau+\gamma\big)
}} \end{eqnarray*}
 with \begin{eqnarray*}w(\tau):=\frac{\varrho^2-1}{(\tau+\varrho)^2}\end{eqnarray*}
 and $\alpha$, $\beta$ and $\gamma$ as in~\eqref{alphabetagamma}. \\
$a_1)$ Let $\varrho\in {\mathbb R}$, with $|\varrho|>1$.
 As, in such case,  $s_\varrho=[-1,1]$, we immediately have  the thesis.
  \\
$a_2)$ Let $\varrho\in {\mathbb R}$, with $|\varrho|<1$. 
We show that $\alpha$ and $\beta$ do not simultaneously vanish. By 
\eqref{conditions5}, $P(z)$ does not have double roots in $(-1,1)$.  If $\alpha=\beta=0$, then 
$b^2-ac=\beta^2-\alpha\gamma=0$ and  $\varrho\in (-1,1)$ satisfies \begin{eqnarray*}a\varrho^2-2b \varrho+ c=\alpha(\varrho^2-1)=0\end{eqnarray*}
which contradicts that $P$ does not have double roots in $(-1,1)$. As $\alpha$ and $\beta$ do not simultaneously vanish, and $s_\varrho=(-\infty, -1]^-\cup[1, +\infty)^-$,  we get
\begin{eqnarray*}
 {\mathcal I}_1&=&\int_{(-\infty, -1]^-\cup [1, +\infty)^- }\frac{d\tau}{\sqrt{(1-\tau^2)
 \big(\alpha\tau^2-2\beta\tau+\gamma\big)
 }}
 \end{eqnarray*}
 As the transformation $\tau=\Gamma_\varrho (z)$ sends roots of $P$ to roots of $R$ and (due to $\varrho\in (-1,1)$) the half--planes $\{\Im z>0\}$, $\{\Im z<0\}$ to $\{\Im \tau<0\}$, $\{\Im \tau>0\}$, respectively, then at least one of the half--planes $\{\Im \tau>0\}$, $\{\Im \tau<0\}$ is free from roots of
$R$. Then we can use  Cauchy theorem, choosing a path integral made of the real axis and half--circle centered at $\tau=0$ with radius $r\to \infty$ on the side of the complex plane where no roots of $R$ occur. Then \begin{eqnarray*}
\int_{(-\infty, -1]^-\cup [1, +\infty)^- }\frac{d\tau}{\sqrt{(1-\tau^2)
 \big(\alpha\tau^2-2\beta\tau+\gamma\big)
 }}=\int_{[-1,1]^+ }\frac{d\tau}{\sqrt{(1-\tau^2)
 \big(\alpha\tau^2-2\beta\tau+\gamma\big) }}
 \end{eqnarray*} 
which concludes.\\
$a_3)$ Let $\varrho\in {\mathbb C}\setminus{\mathbb R}$ and assume that no root of $P$ lies in  ${\mathcal A}(\varrho)$. Then  no root of $R$ lies in  ${\mathcal A}(-\varrho)$. As $s_\varrho={\mathcal S}(-\varrho)$  has as end-points  $-1$ and $1$, by the absence of roots of $R$ in ${\mathcal A}(-\varrho)$, we can deform the path integral, from ${\mathcal S}(-\varrho)^+$ to $[-1,1]^+$.
 $\quad \square$
 \subsection{Proof of ~\eqref{second equality}} In this section we assume ${\rm M}\ge0$. When $\rm M=0$, the part of proof below involving the set ${\mathbb P}^{(3)}_0=\emptyset$ is to be disregarded. The proof is organized in two steps.
\paragraph{\it Step 1}  We prove that
\begin{eqnarray}\label{Step1} \tau_{\rm M}=C
\left\{
\begin{array}{lll}
\displaystyle\int_{-1}^{t_0}\frac{dt}{\sqrt{a_2(1-t^2)(t_0-t)}}\quad &{\rm if}\ ({\rm d}, {\rm f})\in{\mathbb P}^{(3)}_{\rm M}\cup {\mathbb P}^{(4)}_{\rm M}\\\\
\displaystyle\frac{1}{2}\int_{-1}^{1}\frac{dt}{\sqrt{a_2(1-t^2)(t_0-t)}}&{\rm if}\ ({\rm d}, {\rm f})\in{\mathbb P}^{(2)}_{\rm M}
\end{array}
\right.\end{eqnarray}
with
\begin{eqnarray*}C:=2\sqrt{\frac{2\rm d}{v_0}}\,,\quad a_2:=2\sqrt{{\rm f}^2-{\rm M}^2{\rm d}^2}\,,\quad t_0:=\frac{{\rm d}^2-2{\rm f}}{2\sqrt{{\rm f}^2-{\rm M}^2{\rm d}^2}}\,.\end{eqnarray*}
We denote as \begin{eqnarray*}a_1:=2\sqrt{{\rm f}^2-{\rm M}^2{\rm d}^2}+{\rm d}^2-2{\rm f}\,,\quad c_1:=2\sqrt{{\rm f}^2-{\rm M}^2{\rm d}^2}-{\rm d}^2+2{\rm f}\,.\end{eqnarray*}
Consider the triple of polynomials
\begin{eqnarray}\label{P1Q1R1}
P_1(z)&=&a_1z^2+c_1\nonumber\\
Q_1(z)&=&2\left(\sqrt{{\rm f}^2-{\rm M}^2{\rm d}^2}-{\rm f}\right)z^2+2\left(\sqrt{{\rm f}^2-{\rm M}^2{\rm d}^2}+{\rm f}\right)\nonumber\\
R_1(z)&=&{\rm d}^2z^2-4{\rm d}{{\rm M}}z+4{\rm f}-{\rm d}^2
\end{eqnarray}
which verify~\eqref{conditions2} to~\eqref{conditions4}. The roots $\varrho_{1\pm}$ of $Q_1$ are real, opposite and verify $-\varrho_{1-}=\varrho_{1+}>1$ for all $({\rm d}, {\rm f})\in {\mathbb P}^{(2)}_{\rm M}\cup {\mathbb P}^{(3)}_{\rm M}\cup {\mathbb P}^{(4)}_{\rm M}$. Then, for all $x_0\in {\mathbb R}$  and for all $({\rm d}, {\rm f})\in {\mathbb P}^{(2)}_{\rm M}\cup {\mathbb P}^{(3)}_{\rm M}\cup {\mathbb P}^{(4)}_{\rm M}$,
\begin{eqnarray*}
\int_{x_0}^1\frac{dz}{\sqrt{(1-z^2)\left(a_1z^2+c_1\right)}}=
\int_{y_0}^1\frac{dz}{\sqrt{(1-z^2)\big({\rm d}^2z^2-4{\rm d}{{\rm M}}z+4{\rm f}-{\rm d}^2\big)}}
\end{eqnarray*}
 where  $y_0=\Gamma_{\varrho_{1+}(x_0)}$. We next distinguish three cases. In all cases, we denote as $-z_{1-}=z_{1+}=\sqrt{-\frac{c_1}{a_1}}$ the roots of $P_1$, and we use the change of coordinates (between real variables)
\begin{eqnarray}\label{change*****}x=\sqrt{\frac{t-1}{1+t}}:\ [1, +\infty)\to [0, 1)\,.\end{eqnarray} 
\begin{itemize}
\item[\tiny\textbullet] $({\rm d}, {\rm f})\in{\mathbb P}^{(2)}_{\rm M}$.  In this case, \begin{eqnarray*}a_1>0\,,\ c_1<0\,,\ t_0\in (1, +\infty)\,,\ -z_-=z_+\in (0, 1)\,.\end{eqnarray*} 
We take $x_0=z_{1+}$ the maximum root of $P_1$, so $y_0=y_{1+}=\frac{2}{\rm d}\left({\rm M}+\sqrt{{\rm M}^2+\frac{{\rm d}^2}{4}-{\rm f}}\right)$ is the maximum root of $R_1$. We have the equality
 \begin{eqnarray*}
\tau_{\rm M}&=&C\int_{y_{1+}}^1\frac{dz}{\sqrt{(1-z^2)\big({\rm d}^2z^2-4{\rm d}{{\rm M}}z+4{\rm f}-{\rm d}^2\big)}}\nonumber\\
&=&C\int_{z_{1+}}^1\frac{dz}{\sqrt{(1-z^2)\left(a_1z^2+c_1\right)}}\end{eqnarray*}
where the last integral is meant in real sense. 
The change~\eqref{change*****}
transforms the right hand side into
\begin{eqnarray*}
\tau_{\rm M}=
\frac{C}{2}\int_{t_0}^{+\infty}\frac{dt}{\sqrt{a_2(t^2-1)\left(t-t_0\right)}}=\frac{C}{2}\int_{-1}^{1}\frac{dt}{\sqrt{a_2(1-t^2)\left(t_0-t\right)}}
\end{eqnarray*}
where we have used a modification of the integral path, after using that $t_0\in (1, +\infty)$ while $({\rm d}, {\rm f})\in {\mathbb P}^{(2)}_{\rm M}$ and that the function under the integral has only real branching points  and $(-1,1)$, $(t_0,+\infty)$ are the only intervals in the real axis where the function under the integral is real--valued.
\item[\tiny\textbullet] $({\rm d}, {\rm f})\in{\mathbb P}^{(3)}_{\rm M}$. In this case, \begin{eqnarray*}a_1<0\,,\ c_1>0\,,\ t_0\in (-1,0]\,,\ -z_{1-}=z_{1+}\in (1, +\infty)\,.\end{eqnarray*} We take $x_0=-1$, so $y_0=-1$. 
We have the equality
 \begin{eqnarray*}
\tau_{\rm M}&=&C\int_{-1}^1\frac{dz}{\sqrt{(1-z^2)\big({\rm d}^2z^2-4{\rm d}{{\rm M}}z+4{\rm f}-{\rm d}^2\big)}}\nonumber\\
&=&C\int_{-1}^1\frac{dz}{\sqrt{(1-z^2)\left(a_1z^2+c_1\right)}}\nonumber\\
&=&2C\int_{0}^{1}\frac{dx}{\sqrt{(1-x^2)\left(a_1x^2+c_1\right)}}
\end{eqnarray*}
where the last integral is meant in real sense. In the last step, we have used $-z_{1-}=z_{1+}\in (1, +\infty)$.
With the change~\eqref{change*****}
the right hand side becomes
\begin{eqnarray*}
C\int_{1}^{+\infty}\frac{dt}{\sqrt{a_2(t^2-1)\left(t-t_0\right)}}&=&C\int_{-1}^{t_0}\frac{dt}{\sqrt{a_2(1-t^2)\left(t_0-t\right)}}\nonumber\\
\end{eqnarray*}
Here, we have: simplified the factor $2$;  interpreted the integral as a complex one and modified the integral path, after using that $t_0\in (-1,0]$ while $({\rm d}, {\rm f})\in {\mathbb P}^{(3)}_{\rm M}$ and that the function under the integral has only real branching points  and $[1,\infty)$, $[-1,t_0]$  are the only intervals in the real axis where the function under the integral is real--valued.

\item[\tiny\textbullet] $({\rm d}, {\rm f})\in{\mathbb P}^{(4)}_{\rm M}$.  In this case,\begin{eqnarray*}a_1>0\,,\ c_1>0\,,\ t_0\in (-1,1)\,,\ -z_-=z_+\in \ii{\mathbb R}\end{eqnarray*} We take $x_0=-1$, so $y_0=-1$. We have the equality
 \begin{eqnarray*}
\tau_{\rm M}&=&C\int_{-1}^1\frac{dz}{\sqrt{(1-z^2)\big({\rm d}^2z^2-4{\rm d}{{\rm M}}z+4{\rm f}-{\rm d}^2\big)}}=C\int_{-1}^1\frac{dz}{\sqrt{(1-z^2)\left(a_1z^2+c_1\right)}}\nonumber\\
&=&2C\int_{0}^1\frac{dz}{\sqrt{(1-z^2)\left(a_1z^2+c_1\right)}}\end{eqnarray*}
where the last integral is meant in real sense. 
In the last step, we have used $-z_{1-}=z_{1+}\in \ii{\mathbb R}$.
The change~\eqref{change*****}
transforms the right hand side into
\begin{eqnarray*}
\tau_{\rm M}=C\int_{1}^{+\infty}\frac{dt}{\sqrt{a_2(t^2-1)\left(t-t_0\right)}}=C\int_{-1}^{t_0}\frac{dt}{\sqrt{a_2(1-t^2)\left(t_0-t\right)}}\end{eqnarray*}
again by a modification of the integral path, after using  $t_0\in (-1, 1)$ while $({\rm d}, {\rm f})\in {\mathbb P}^{(4)}_{\rm M}$ and that the function under the integral has only real branching points  and $(-1,t_0)$, $(1,+\infty)$ are the only intervals in the real axis where the function under the integral is real--valued.

\end{itemize}
Then~\eqref{Step1} is completely proved.
\paragraph{\it Step 2} We now consider the polynomials
\begin{eqnarray*}
P_2(t)&=&Q_2(t)={\rm M}^2t^2-2({\rm M}^2-{\rm f})t+{\rm M}^2-2{\rm f}+{\rm d}^2\nonumber\\
R_2(t)&=&a_2(t_0-t)=-2t\sqrt{{\rm f}^2-{\rm M}^2{\rm d}^2}-2{\rm f}+{\rm d}^2\nonumber\\
\end{eqnarray*}
which verify~\eqref{conditions2} to~\eqref{conditions4}.
As  $P_2=Q_2$ has real--valued roots 
while $({\rm d}, {\rm f})\in {\mathbb P}^{(2)}_{\rm M}\cup {\mathbb P}^{(3)}_{\rm M}\cup {\mathbb P}^{(4)}_{\rm M}$, we have 
\begin{eqnarray*}
\int_{-1}^{t_0}\frac{dt}{\sqrt{(1-t^2)\left(-2t\sqrt{{\rm f}^2-{\rm M}^2{\rm d}^2}-2{\rm f}+{\rm d}^2\right)}}&=&\int_{-1}^{z_0}\frac{dz}{\sqrt{(1-z^2)\left({\rm M}^2z^2-2({\rm M}^2-{\rm f})z+{\rm M}^2-2{\rm f}+{\rm d}^2\right)}}\nonumber\\
&=&\int_{0}^{w_0}\frac{d w}{\sqrt{w(2-w)\left({\rm M}^2w^2-2{\rm f}w+{\rm d}^2\right)}}
\end{eqnarray*}
having let $w:=1-z$, and having named and $z_0$ the minimum root of ${\rm M}^2z^2-2({\rm M}^2-{\rm f})z+{\rm M}^2-2{\rm f}+{\rm d}^2$ and
$w_0$ the minimum root of ${\rm M}^2w^2-2{\rm f}w+{\rm d}^2$. $\quad \square$

\begin{remark}\rm
One might wonder if the method of proof of~\eqref{second equality} would work also in the case $({\rm d}, {\rm f})\in {\mathbb P}^{(1)}_{\rm M}$, which has been treated separately in Section~\ref{Proof of Theorem1.2}, via
 dynamical arguments. As a matter of fact, when $({\rm d}, {\rm f})\in {\mathbb P}^{(1)}_{\rm M}$, ${\mathcal C}_{\varrho_{1\pm}}$ is the unit circle centered at $0$ and the roots of $P_1$ are complex and both belong to ${\mathcal C}_{\varrho_{1\pm}}$. In that case, there is no room to apply Proposition~\ref{Legendre polynomials2COR}, with the  choice of $P_1$, $Q_1$, $R_1$ in~\eqref{P1Q1R1}.
\end{remark}
 \appendix
 
\subsubsection*{Acknowledgments}
I am grateful to G. Tommei for useful talks, and to the anonymous reviewers for their insights, which considerably  helped me to improve the readability of the paper. This work is part of the ERC project  Stable and Chaotic Motions in the Planetary Problem, 2016--2022 (Grant 677793).
\addcontentsline{toc}{section}{References}
 \bibliographystyle{plain}

\end{document}